\documentclass[11pt]{amsart}
\usepackage[all]{xy}

\newtheorem{definition}{Definition}
\newtheorem{theorem}{Theorem}

\newcommand{\he}{\simeq}                %homotopy equivalent

\newcommand{\x}{\mathcal{F}_{l}}
\newtheorem{lemma}{Lemma}[section]
\newtheorem{corollary}[lemma]{Corollary}
\newtheorem{proposition}[lemma]{Proposition}
\theoremstyle{remark}
\newtheorem{remark}[lemma]{Remark}
\newtheorem{example}{Example}
\newcommand{\bib}[6]{ \bibitem{#1}  #2, \textit{ #3\/}, #4 \textbf{#5} #6}

\begin{document}

\title[Hopf invariants and hamiltonians]{Homotopical
Dynamics IV: Hopf invariants and hamiltonian flows}
\author{Octavian Cornea}
\address{Universit\'{e} de Lille 1\\
U.F.R. de Math\'{e}matiques \& U.M.R. 8524\\
59655 Villeneuve D'Ascq, France}

\email{cornea@gat.univ-lille1.fr}
\urladdr{http://www-gat.univ-lille1.fr/\~{}cornea/octav.html}

\date{22.04.2001}

\maketitle

\section*{Introduction.}

Let $(M,\omega)$ be an $n$-dimensional, symplectic,  not necessarily compact
manifold. Suppose that $f$ is a $C^{2}$, real  function on $M$. A fundamental question
 is whether the hamiltonian flow $h^{f}$ induced by $(f,\omega)$ has periodic orbits.
Many results are known when $M$ or some of the level surfaces of $f$ have compact
components or, at least, are of finite volume.
 However, if this is not the case - and our focus will be
on this non-compact case - very few results exist ( see \cite{Ho} \cite{HoZe} as well as the
discussion at the beginning of \S\ref{sec:orb}).

\

In this non-compact context the first natural step in the search for periodic orbits
is to detect {\em bounded} ones. The key underlying idea of this paper is that
 algebraic topological constraints imposed to certain
 invariants associated to a {\em gradient} flow of $f$ imply the existence of bounded orbits for
$h^f$. Once the existence of bounded orbits is established,
under favorable circumstances, application of the $C^{1}$-closing lemma
leads to  periodic ones.

\

Fix a riemannian metric $\alpha$ on $M$.  For technical reasons, we will assume in this
paper that $M$ {\em is simply-connected}. Let $S$ be an isolated, compact, invariant set of the
negative gradient flow $\gamma(f)$ of $f$ \cite{Con}. For such an $S$  we consider two
associated invariants.  The first, $c(S)$, is the Conley index \cite{Con}. The second invariant
$d(S)\in {\bf Z}/2$  is new and of a different type. Here are some consequences of our results:

{\em  Assume that the negative gradient flow of the function $f:M\to R$ has an isolated
invariant set $S$ such that $d(S)=1$ and $S$ satisfies a technical property
$(\ast )$ (that will be discussed below).

\begin{itemize}
\item[(i)]  If $f$ is Morse, then there are infinitely many regular hypersurfaces $f^{-1}(a)$
which after deformation by an arbitrarily small, compactly supported isotopy carry a closed
characteristic.
\item[(ii)] For a general $f$ there is some $C^{2}$ -neighbourhood of
$f$ in the Whitney (strong) topology containing a dense family of functions each of whose
hamiltonian flow has infinitely many distinct, periodic orbits.
\end{itemize} }

These results are non-trivial because, in general, the hypersurfaces $f^{-1}(a)$ appearing above
have no compact components.  Some explicit such examples of
connected, non-compact  hypersurfaces $V$  are given in  \S\ref{subsec:appl}. We shall also
see there that there are $(f,S)$  such that property $(\ast)$ is satisfied, $h^{f}$  has no
non-trivial {\em bounded} orbit and the family of hamiltonians that have some non-trivial
bounded orbit is not dense in any neighbourhood of $f$ (obviously in this case $d(S)=0$).
In fact, under the assumptions at (i) above we obtain that all hypersurfaces $f^{-1}(a)$
that appear in that statement (before deformation) contain at least one bounded orbit
of $h^{f}$.

In many cases, $d(S)$  can be computed homotopically
out of an index pair $(N_{1},N_{0})$ (in the sense of the Conley index theory, see
\cite{Con} \cite{Sa} and also  \S\ref{subsubsec:conlind}) of $S$. By definition, $c(S)$ is
the homotopy type of $N_{1}/N_{0}$. As for $d(S)$, an important property is that it equals
$1$ if one of the following two conditions is satisfied: a. $c(S)$ has torsion in homology; b.
there exists an index pair for $S$ with $N_{1}$ and $N_{0}$  both simply-connected, the
homology of $c(S)$ is torsion free and a certain chain complex
$\mathcal{C}_{min}(N_{1},N_{0})$ of free
$\pi_{\ast}^{S}((\Omega M)^{+})$-modules canonically associated to a minimal $CW$-cell
decomposition of the pair $(N_{1},N_{0})$ has a non trivial differential.
We call this complex, which is unique up to isomorphism and is closely related to the twisted
complexes introduced by Baues \cite{Ba} in homotopy theory, the minimal Hopf complex of
$(N_{1},N_{0})$. This name is motivated by the fact that the differential is expressed
in terms of certain relative Hopf invariants associated to  succesive cell-attachments in a
minimal cell-decomposition of $(N_{1},N_{0})$. There are reasonably simple criteria
implying the non-vanishing of the differential, in particular one in terms of the non-vanishing
of certain differentials in the Atiyah-Hirzebruch-Serre $\pi_{\ast}^{S}(\Omega
M^{+})$-spectral sequence of the fibration induced over $(N_{1},N_{0})$ from the path-loop
fibration $\Omega M\to PM\to M$.

Here is now property $(\ast )$: we say that $S$ satisfies the property $(\ast )$  (with respect
to $v\in {\bf R}$) if for a value $v\in{\bf R}$  there exists an index pair
$(N_{1},N_{0})$ for $S$ (in the flow $\gamma(f)$) such that $N_{0}\subset f^{-1}(v)$
and $f$ is regular on $N_{0}$.
There are many simple criteria that insure that this condition is satisfied (see Lemma
\ref{lem:prop_ast}).

\

We explain below the ideas used in the proof of the results mentioned above and, at the
same time, describe the contents of the paper.

Given a gradient flow of $f$ and an isolated invariant set $S$ in this flow, Conley's theory
provides a way to associate some homotopy theoretical meaningful notions to $(f,S)$.
In particular, because we are talking here about a gradient flow
on a finite dimensional manifold, isolating neighbourhoods as well as index pairs of $S$ are
generally not difficult to determine. Assuming such an index pair $(N_{1},N_{0})$ constructed
the main aim becomes to translate the homotopical properties of this pair into recurrence
properties of $h^{f}$.

Notice that
even if the function $f$  is not Morse itself it can be morsified in a neighbourhood of $S$. In
other  words,  given an isolating neighbourhood $N$ of $S$  we may continue $S$ inside $N$
to an isolated invariant set $S'$ of the gradient flow of a function $g$ which agrees with $f$
outside $N$ and satisfies the Morse-Smale condition inside $N$. Such a $g$ can be constructed
to be as close to $f$ in the $C^{2}$-strong topology as desired.
Consider couples $(g,N)$ as above. Explicitely, $N$
is an isolating neighbourhood of the negative gradient flow of $g$ and, inside $N$, $g$ is
Morse and satisfies the Morse-Smale transversality condition.
To such a couple - which will be called a {\em local Morse function} - we associate
certain invariants which make the connection between the homotopy type
of the pair $(N_{1},N_{0})$ and the existence of bounded orbits of $h^{g}$.

These invariants are Morse-type complexes only that,
 instead of measuring the space of flow lines joining
critical points (of $g|_{N}$)
by just counting the number of its elements (when this number is finite), we
consider consecutive critical points of indexes different by possibly more than one but
successive in a fixed index set $I\subset {\bf N}$ and measure the corresponding
space of  connecting flow lines (called connecting manifold) by a certain framed bordims class
that belongs to $\Omega_{\ast}^{fr}(\Omega M)$ (the precise definition is in
\S\ref{subsec:rec}). We denote the resulting $\Omega_{\ast}^{fr}(\Omega M)$-chain complexes
 by $C^{I}_{N}(g)$ and we call them
{\em extended Morse complexes} of $g$ (relative to $N$).
The invariant $d(S)$ is defined in terms of these complexes
 (see \S \ref{subsec:appl}): $d(S)=1$ if there is
a $C^{2}$ -neighbourhood, $\mathcal{U}_{f}$, of $f$ such that
for each  morsification $(g,N)$ of $S$  with $g\in\mathcal{U}_{f}$
there is an extended Morse complex $C^{I}_{N}(g)$  with
non-trivial differential.

\

The paper is structured by the point of view that extended Morse complexes,
introduced in Section \ref{subs:def_compl},
provide a bridge between homotopy theory and recurrence properties of
the respective hamiltonian flows. The two other
sections of the paper deal each with one of the two sides of the story.

\

The relation between extended Morse complexes and
homotopy theory is the task of
Section \ref{sec:conn}. Assume $(g,N)$ is a local Morse function and let $P$ and $Q$ be
consecutive critical points of $g|_{N}$  (which means $p=ind(P)>ind(Q)=q$ and $N$
does not contain broken flow lines of the gradient flow of $g$ joining $P$ to $Q$). The
geometry of the connecting manifold of $P$ and $Q$ produces, by the Thom-Pontryaguin
construction, (see again \S \ref{subsec:rec}) a homotopy class  $h(P,Q)\in
\pi_{p-1}(\Sigma^{q}(\Omega M)^{+}) )$ whose stabilization is the bordism class mentioned
above.  The main result of this section shows that $h(P,Q)$ equals a certain homotopical object
called the {\em relative Hopf invariant} associated to the cell attachments
corresponding to the ``passage" through
$P$ and $Q$ (this extends results in \cite{Fr} and  \cite{Co2}).  In particular, this  allows the
identification of the extended Morse complexes of $(g,N)$ with Hopf - type
complexes associated to a relative, not necessarily minimal,
cell-decomposition of the index pair $(N_{1},N_{0})$ (recall that
$\Omega_{\ast}^{fr}(A)\approx \pi_{\ast}^{S}(A^{+})$; $A^{+}=(A)^{+}=A$ union a
disjoint point).
As a corollary of this identification we prove in
\S\ref{sec:rig} some rigidity properties of these extended Morse complexes. In particular,
returning  to the function $f$ and the isolated invariant set $S$, if the index pair
$(N_{1},N_{0})$ has torsion free integral relative homology and both $N_{1}$ and
$N_{0}$ are simply connected, then
$\mathcal{C}_{min}(N_{1},N_{0})$ is defined and, if the
differential in $\mathcal{C}_{min}(N_{1},N_{0})$  is non-trivial, then the differential of  the
complex $C^{I_{g}}_{N}(g)$ does not vanish either for any morsification $(g,N)$ of $S$
and  some index set $I_{g}$. It also follows that if $c(S)$ has torsion in homology, then the
differential of such a complex  $C^{I_{g}}_{N}(g)$ is again non-trivial. This leads to the
fact that in both these two cases $d(S)=1$.

\

In Section \ref{sec:orb} we relate the differential of extended Morse complexes
 to hamiltonian flows. It is only here that property $(\ast)$ is of use. For a local Morse function $(g,N)$
whose maximal invariant set $S'$ inside $N$ satisfies condition $(\ast)$ the main result of this
section shows that
the non-vanishing of the differential in $C^{I_{g}}_{N}(g)$  implies that the hamiltonian flow
of $g$ has some bounded orbits  situated on regular hypersurfaces.
The (very) rough idea of the proof is as follows. The non-vanishing of the differential above implies
that some bordism class [$h(P,Q)$] is non-trivial. Let $V=g^{-1}(a)$ be such that  it separates
$P$ and $Q$, $g(P)>g(Q)$. We show that if the flow $h^{g}$ has no bounded
orbits on $V$, then the connection between $P$ and $Q$ given by the negative gradient
flow lines of $g$ can be annihilated by letting the unstable sphere of $P$ in
$V$, $W^{u}(P)\cap V$, be deformed by (a perturbation of) $h^{g}$ till it
becomes disjoint from the stable sphere of $Q$ in $V$, $W^{s}(Q)\cap V$. This leads to
a contradiction because it would imply that the bordism class $[h(P,Q)]$ is trivial.
It is useful to note that this type of argument actually does apply to any flow orthogonal to the
gradient flow of $g$. We use the $C^{1}$-closing lemma of Pugh and Robinson \cite{PuRo} to
deduce the existence of periodic orbits for the flow $h^{g'}$ for some functions $g'$ close
to $g$ as well as to obtain the existence of closed characteristics on some hypersurfaces close to
those of $g$. This immediately leads to a proof of the points (i) and (ii) above.
When $d(S)=0$ or the condition $(\ast)$ is not satisfied a different method is available. This
continues an idea that first appeared  in \cite{Co3}.  It makes use of the dynamical
interpretation of Spanier-Whitehead duality first discussed in \cite{Co1} to show that, under
certain circumstances, the lack of self-duality of  the Conley index of a compact invariant set in
a gradient flow can force the appearance of some non-constant, bounded orbits for the
associated hamiltonian flow.
A result obtained in this way claims the existence of a family of functions dense in a
neighbourhood of $f$ each of whose induced hamiltonian flow has at least as  many closed
(possibly homoclinic), distinct, non-trivial orbits as $\sum_{k\not=n/2}rk (H_{k}(c(S);{\bf
Z})$.

\

Section \ref{sec:conn} applies to local Morse functions on a manifold $M$ that is
not necessarily  symplectic. The results contained here are of interest independentely of their
applications in Section \ref{sec:orb}.  Indeed, at the center of this series of papers is the old
idea, going back to the work of Poincar\'e, that homotopy theory is relevant in
understanding flows. The first systematic attempt to go beyond applications of homology in this
direction have been provided by Conley's theory \cite{Con} and by work of John Franks,
particularly in \cite{Fr}. A key point is that, essentially, Conley's approach only provides
information about the simplest existing flow inside an isolating neighbourhood if the boundary
behaviour is fixed. This means that it is not only natural but also quite meaningful to first try to
understand the simplest isolated invariants sets possible - isolated invariant points in gradient
flows - and the simplest attractor-repellor pairs possible - which are consecutive critical points
for Morse-Smale functions together with the gradient flow lines that connect them. Isolated
singularities of real functions have been studied from this perspective in \cite{Co3}.
In \cite{Fr} John Franks  has started the application of more refined homotopical tools to the
study of  these ``simplest" attractor-repellor pairs formed by consecutive critical points in
Morse-Smale gradient flows. His results have been extended
by the author in \cite{Co2} and this work is continued and strengthened here.
 Extended Morse complexes provide a tool to encode the data coming from these
attractor-repellor pairs in an algebraic way and, therefore, their study seems to be a
worthwhile topic in itself. Because of this we include
 in Section \ref{sec:conn} more properties of these complexes than those strictly needed
later in the paper. In particular, we discuss a connection between this topic
and the problem of the existence of smoothings of Poincar\'e duality spaces (along the lines of
\cite{Co2}).

\

Most of the paper is written in the language of Conley's index theory.
A notable exception is  \S \ref{subsec:hopf} which concentrates all the non-trivial homotopy
theoretical arguments of the paper. This sub-section is of some intrinseque interest and can be
read independently of the rest of the paper.

To make
the paper as widely accessible as possible all the needed elements of Conley index theory
together with the basics of Morse theory that we use are recalled in Section
\ref{subs:def_compl}.
\tableofcontents

\section{Recalls and preliminary definitions.} \label{subs:def_compl}

\

This section does not contain new results or proofs. It serves two purposes. The first is to fix
notations and standard facts concerning Morse functions and  the Conley index.  The second
is to introduce the definitions of ``generalized" and  ``extended" Morse complexes.

\subsection{The Conley index.}\label{subsubsec:conlind}

\

We  recall some elements of Conley index theory \cite{Con}\cite{Sa}.

Let $\gamma: X\times {\bf R}\longrightarrow X$  be a continuous flow on a locally
compact, metrizable space $X$.  Given a subset $N\subset X$ we denote by
$I_{\gamma}(N)$ or $I(N)$ the maximal invariant set of $\gamma$ included in $N$.
A compact $N\subset X$ is called isolating neighbourhood in the flow $\gamma$
if $I(N)\subset Int(N)$. It follows that $I(N)$ is compact.
In other words a compact $N\subset X$ is an isolating neighbourhood if all orbits of type
$\gamma_{\bf R}(x)$,  $x\in N$ that are included in $N$ do not intersect $\partial N$.
 An invariant set of $\gamma$, $S\subset X$  is isolated if $S=I(N)$ for some isolating
neighbourhood $N$.  Fix now an isolated invariant set $S$ inside some
isolating neighbourhood $N$  in the flow $\gamma$.
A pair $(N_{1},N_{0})$ of
compact subsets of $N$ is an index pair for $S$ in $N$ if $N_{0}\subset N_{1}$,
$N_{1}-N_{0}$ is a  neighbourhood of $S$, $S=I(\overline{Int(N_{1}-N_{0})})$, $N_{0}$
is positively invariant in $N_{1}$  and, if for $x\in N_{1}$ there is some $T\geq 0$ such that
$\gamma_{T}(x)\not\in N_{1}$, then there exists a $\tau >0$,
$\tau<T$ with $\gamma_{t}(x)\in N_{1}$ for $0\leq t\leq \tau$ and $\gamma_{\tau}(x)\in
N_{0}$. For such an index pair $(N_{1},N_{0})$ we will call $N_{0}$ the exit set.
A key basic result is that there are index pairs inside any isolating neighborhood of $S$. The
Conley index of $S$, $c_{\gamma}(S)$,  is the homotopy type of the quotient space
$N_{1}/N_{0}$. It is independent of the choices of the index pair. It is also invariant to
continuation in the sense that if $\gamma^{\lambda}:X\times {\bf R}\longrightarrow X$ is a
family of flows depending  continuously on the parameter $\lambda\in [0,1]$ and
$S$ is an isolated invariant set of $\gamma$ viewed as a flow on $X\times [0,1]$, then
$S_{\lambda}=S\bigcap X\times\{\lambda\}$ is an isolated invariant set of
$\gamma^{\lambda}$ and $c_{\gamma}(S)= c_{\gamma^{\lambda}}(S_{\lambda})$
for all $\lambda\in [0,1]$.
For an arbitrary point $x\in X$ let its $\omega$ and $\omega^{\ast}$ limits be defined
by $$\omega(x)=\{z\in X :\exists \{t_{n}\}\rightarrow +\infty, \
z= \lim_{n\to\infty}\gamma_{t_{n}}(x)\}$$
and $$\omega^{\ast}(x)=\{z\in X : \exists \{t_{n}\}\rightarrow -\infty, \
z=\lim_{n\to\infty}\gamma_{t_{n}}(x)\}$$
If $N$ is an isolating neighbourhood of the isolated invariant
set $S$, then we denote by $$W^{u}_{N}(S)=\{x\in N : \forall  t\in (-\infty,0], \
\gamma_{t}(x)\in N, \ \omega^{\ast}(x)\subset S\}$$
and
$$W^{s}_{N}(S)=\{x\in N : \forall  t\in [0,\infty), \
\gamma_{t}(x)\in N, \ \omega(x)\subset S\}$$

the unstable and stable manifolds of $S$ in $N$.
Both sets are compact. We will generally omit the index
$N$ from the notation if this does not create any ambiguity.
A special type of index pairs that are particularly useful
are regular index pairs \cite{Sa}. Two of their important properties are that $N_{0}$
is a neighbourhood deformation retract and that the arrival time in $N_{0}$ is a continuous
function on $N_{1}-W_{N_{1}}^{s}(S)$ where $S$ is the maximal invariant set inside
$N_{1}$.  We also need the notion of bi-regular index block. This is a triple
$(N_{1},N_{0},N_{2})$ such that
$(N_{1},N_{0})$ and
$(N_{1},N_{2})$ are respectively regular index pairs for the direct and inverse flows and
$N_{0}\cup N_{2}$ includes $\partial N_{1}$  (see \cite{Re}; the inverse flow of $\gamma$
is defined by $-\gamma_{t}(x)=\gamma_{-t}(x)$). In our case, $\gamma$ will always be
defined on a smooth manifold and is differentiable of class at least $C^{1}$.  We will also need
a special type of bi-regular index block that we call {\em strong index block}. This is a triple
$(N_{1},N_{0},N_{2})$ such that $N_{1}$ is a submanifold of the same dimension  as the
underlying manifold, it has a smooth interior, $\partial N_{1}=N_{0}\cup N_{2}$ with
$N_{0}$ and $N_{2}$ smooth manifolds possibly with boundary $\partial N_{0}=\partial
N_{2}=N_{0}\cap N_{2}$  and such that the vector field $X(\gamma)$  associated to
$\gamma$ is never tangent to $\partial N_{1}$ and points strictly inside $N_{1}$ on
$N_{2}$ and strictly outside on $N_{0}$. Such a strong index block can be easily
constructed inside any isolating neighbourhood by making use of Lyapounov functions for the
direct and inverse flows \cite{Sa}. Here is a simple but important property of strong index
blocks. Fix $(N_{1},N_{0},N_{2})$ a strong index block for $\gamma$.  If $\gamma'$ is a
flow such that $X(\gamma')$ is sufficientely $C^{0}$-close to  $X(\gamma)$  on
$\partial N_{1}$, then $(N_{1},N_{0},N_{2})$ is also a strong index block for $\gamma'$.
By using the one parameter family of flows induced $tX(\gamma)+(1-t)X(\gamma')$,
it results that the maximal invariant sets of $\gamma$ and $\gamma'$ inside $N_{1}$ are
related by continuation. This construction is useful in continuing a given isolated invariant set
$S$ of some flow $\gamma$ to simpler invariant sets. In particular, it is easy to show by using
the two Lyapounov functions mentioned above that, on a smooth manifold, any isolated
invariant set can be continued to one  in a gradient flow.

The Conley index also associates some relevant
homotopy data to an attractor-repellor pair $(S;R,A)$ in the flow $\gamma$.
This is an isolated invariant set such that
$A,R\subset S$ are disjoint  isolated invariant sets and any $x\in S-(A\cup R)$  has an
verifies $\omega(x)\subset A$, $\omega^{\ast}(x)\subset R$.  To such an attractor repellor
pair there coresponds a homotopy cofibration sequence \cite{Sa} called the attractor-repellor
cofibration sequence:
$$ c_{\gamma}(A)\to c_{\gamma}(S)\to c_{\gamma}(R)\stackrel{\delta(R,A)}{\to}\Sigma
c_{\gamma}(A)$$
One important feature is that the map $\delta(R,A)$  admits a flow based
 explicit defintion besides the one coming from the Baratt-Puppe sequence. This map is called
the connection map of the attractor repellor pair.

There is a variant of the Conley index that has been recentely introduced by
Mrozek, Reineck and Srzednicki \cite{MrReSr} and which is useful in our
context. For a flow $\gamma$ on $X$ having an isolated invariant set $S$  consider a
regular index pair $(N_{1},N_{0})$ of $S$. Consider the space $X'$ obtained from the
disjoint union of $X$ and  of $N_{1}$ modulo the identification of $N_{0}\hookrightarrow
N_{1}$ with $N_{0}\hookrightarrow X$. Of course $X\subset X'$ and there is also a map
$p:X'\to X$ which is the identity on $X$ and the inclusion on $N_{1}$. The {\em global}
Conley index of $\gamma$ as defined in \cite{MrReSr} is denoted by
$\overline{c}_{\gamma}(S)$ and it is the homotopy type (in the obvious sense) of the
object given by the triple $(X',X,p)$.  The main result in \cite{MrReSr} is that this global
Conley index satisfies all the properties of the usual one, in particular it is independent  of  the
choice of index pairs and is invariant to continuation. All algebraic topological invariants
(i.e. homology or homotopy groups) of $\overline{c}_{\gamma}(S)$ will be
understood to be the respective invariants of the pair $(X',X)$.

The notions above adapt in an obvious way to local flows \cite{Sa}. This is particularly
useful because many of the flows that we consider in this paper are in fact only partially
defined (due to the non-compactness of our underlying manifold) but they
are all local flows.

\subsection{Elements of Morse theory and extended Morse complexes.}\label{subsec:rec}

\

Let $M$ be our fixed, smooth, connected,  not necessarily compact, $n$-dimensional
manifold with a fixed  riemannian metric $\alpha$.
For a $C^{1}$ function $f:M\to {\bf R}$  we denote by $\nabla(f)$  the
$\alpha$ -gradient vector field of $f$ and let $\gamma(f)$ be the (partially defined)
flow induced by $-\nabla(f)$.

\begin{definition}\label{def:loc_m}  A {\em local Morse function} on $M$ is a pair
$(f,N)$ such that $f:M\to {\bf R}$ is smooth, $N$ is an isolating  neighbourhood in the flow
$\gamma(f)$, $f|_{N}$ is Morse (meaning that all critical points  of $f$ inside $N$ have a
non-degenerate hessian) and $f$ verifies the Morse-Smale transversality condition in $N$.
\end{definition}

Here is the meaning of the last assumption.
Let $Crit^{N}(f)$  be the critical points of $f|_{N}$ and let
$Crit_{k}^{N}(f)$ be the set of  those critical points $x$ of index $ind(x)=k$.
For $x\in Crit^{N}_{k}(f)$ the stable and unstable manifolds of $x$ in $N$
$W^{s}_{N}(x)$ and $W^{u}_{N}(x)$ are both defined with respect to the
flow $\gamma(f)$ as in \S\ref{subsubsec:conlind}. Their intersections with the interior of
$N$ are, respectively, open $n-k$ and $k$ -dimensional manifolds.
If $N$ has the property that there are $N_{0}$, $N_{2}$ such that $(N,N_{0},N_{2})$
is a strong index block, then these manifolds are diffeomorphic to open euclidean disks
of the corresponding dimensions.
The Morse-Smale condition requires that for any $x,y\in Crit^{N}(f)$ the intersection of
$W^{u}_{N}(x)$ and $W^{s}_{N}(y)$ be transverse.

It is useful to notice that if $f:M\to {\bf R}$ is $C^{1}$ and $N$ is an isolating
neighbourhood for $\gamma(f)$, then the invariant set $I_{N}(f)=I_{\gamma(f)}(N)$  is
formed by all the critical points in $Crit^{N}(f)$ and all the points situated on some flow
line of $\gamma(f)$ that is contained in $N$ and joins two of these critical points.

In case $(f,N)$ is a local Morse function  the intersection of $W_{N}^{u}(x)\cap
W^{s}_{N}(y)\cap f^{-1}(a)$, $x,y\in Crit^{N}(f)$,
$a\in {\bf R}$ regular value of $f$, is void or a compact stratified manifold with the top
dimensional stratum of dimension $ind(x)-ind(y)-1$.
Two critical points $P,Q\in Crit^{N}(f)$  are called consecutive (in $N$) if $f(P)>f(Q)$ and if
$\overline{W^{u}_{N}(P)}\bigcap\overline{W^{s}_{N}(Q)}\bigcap
Crit^{N}(f)-\{P,Q\}=\emptyset$. In this case the intersection $Z_{N}(P,Q)=
W_{N}^{u}(x)\cap W^{s}_{N}(y)\cap f^{-1}(a)$ is a compact manifold of dimension
$ind(P)-ind(Q)-1$ and its diffeomorphism type does not depend of $a$. We will use the
notation $Z(P,Q)$ for this manifold if $N$ is fixed.

Extended Morse complexes are a convenient way to encode some of the characteristics
of the  connecting maniflods of consecutive critical points.
A very simple, general definition is the following.

\begin{definition}\label{morsecomplex} Given
 some graded ring $\mathcal{R}_{\ast}$, a generalized Morse complex
associated to a local Morse function $(f,N)$ is a chain complex of free
$\mathcal{R}_{\ast}$-modules
$(\mathcal{C},d)$ such that $\mathcal{C}=\mathcal{R}_{\ast}\otimes {\bf
Z}[Crit^{N}_{\odot}(f)]$  (where
$Z[A_{\odot}]$ is the graded, free abelian group generated by $A_{\odot}$) and for $x\in
Crit^{N}_{i}(f)$, $$dx=\sum_{y\in Crit^{N}_{<i}(f)}d(x,y)y$$ $d(x,y)=0$ if $x$ and $y$
are not consecutive.
\end{definition}

The definition above becomes interesting whenever the coefficients $d(x,y)$
reflect part of the geometry of the flow $\gamma(f)$. Typically,
these coefficients "measure" the connecting manifolds $Z(x,y)$.
At the same time, richer the ring $\mathcal{R}$ is, finer is the structure encoded in the
complex.

The classical example of a Morse complex due to Thom, Milnor,
Smale appears when $\mathcal{R}={\bf Z}$ and $d(x,y)$ is null except when
$ind_{f}(x)=ind_{f}(y)+1$ and in that case $d(x,y)$ is equal to the number $\#(x,y)$
(or $\#_{f}(x,y)$)  of points in $Z(x,y)$ counted with appropriate signs.
We will denote this complex by $C^{Mo}_{N}(f)$.
Another  related example that has received significant attention recentely is the
Morse-Novikov complex \cite{No}.

The key example that interests us here is an extension of the classical Morse case.
 Consider the local Morse function $(f,N)$ and let $P\in Crit^{N}(f)$. Let $S^{u}(P)$
(respectively $S^{s}(P)$) be the unstable sphere (resp. stable sphere) of
$P$ which is defined by $S^{u}(P)=W^{u}_{N}(P)\bigcap f^{-1}(a)$ (resp.
$S^{s}(P)=W^{s}_{N}(P)\bigcap f^{-1}(a)$) for $a$ in between $f(P)$ and
$f(P)+\epsilon$ with $\epsilon$ small in module and negative (respectively positive).
 The Morse lemma shows that if $\epsilon$ is sufficientely small, then
$S^{u}(P)$ and $S^{s}(P)$ are both diffeomorphic to spheres of dimensions
respectively $ind(P)-1$ and $n-ind(P)-1$. Assume now that $Q\in Crit^{N}(f)$
such that $P,Q$ are consecutive. By transporting $Z(P,Q)$ along the flow
$-\gamma(f)$ we see that $Z(P,Q)$ is embedded in $S^{u}(P)$.
It is well known \cite{Fr} that $Z(P,Q)$ is actually framed in
$S^{u}(P)$ with a standard normal framing.  There is also a map $$j(P,Q):Z(P,Q)\to\Omega
M$$ defined as follows (here, and in the rest of the paper, $\Omega M$ is the space of
{\em pointed loops on $M$}). First fix  a path $w$ in $M$ joining $Q$ to $P$.  The map $j(P,Q)$
associates to each point $z\in Z(P,Q)$ the loop obtained by following from $P$ to $Q$ the
flow line of $\gamma(f)$  that passes through $z$ and returning to $P$
via $w$. The homotopy type of this map is well defined and independent of
the choice of $w$ up to conjugation by elements in $\pi_{1}(M)$.
To avoid this ambiguity as well as further technical complications, we assume in this paper, as
mentioned in the introduction, that $M$ {\em is simply-connected}. In this case the framed
bordism classes $[Z(P,Q)]\in\Omega^{fr}_{p-q-1}(\Omega M)$ of the connecting manifolds
$Z(P,Q)$  are well defined. By convention, we will put $[Z(x,y)]=0$ whenever $x$ and $y$
are not consecutive critical points.

\begin{definition}\label{extendedmorse} For a local Morse function $(f,N)$  let
$I=(a_{1},\ldots a_{i}\ldots )\subset {\bf N}$ be a strictly increasing sequence such that if
$Crit_{k}^{N}(f)\not=\emptyset$, then $k\in I$. The extended Morse complex of $(f,N)$
relative to $I$ is a generalized Morse complex $(C^{I}_{N}(f),d^{I}_{f})$  with
$\mathcal{R}_{\ast}=\Omega^{fr}_{\ast}(\Omega M)$ and $d^{I}_{f}(x, y)=[Z(x,y)]$  if for
some $i$, $ind_{f}(y)=a_{i}$, $ind_{f}(x)=a_{i+1}$  and
$d^{I}_{f}(x,y)=0$ otherwise.
\end{definition}

To understand this object, the key point is to relate $d^{I}_{f}$ to homotopy theory.
This is  the purpose of the main theorem in \cite{Co2} (in the case when $N=M$ which
corresponds to the standard case of a Morse-Smale function on a compact manifold)
and the results obtained there will be strengthened considerably here in Section 2 where
we also show $(d^{I}_{f})^{2}=0$.

\begin{remark}\label{rem:incompl} a. Notice that not all bordism classes
of connecting manifolds of consecutive critical points are contained in a Morse complex of type
$(C^{I}_{N}(f),d^{I}_{f})$ even if the index set $I$ is allowed to vary. This happens
because two critical points $P$ and $Q$ might be consecutive  with respect to the flow
$\gamma(f)$  even if there might be a third critical point $R\in N$ of $f$ such that
$ind(Q)<ind(R)<ind(P)$.

b. The index sets $I$ above are useful in the definition of
extended Morse complexes because the complexes associated to two
different functions that are related by continuation and constructed with respect to
the same index set $I$ will be shown to be comparable by a chain map in \S\ref{sec:rig}.

c. There is an obvious ring map $r:\Omega^{fr}_{\ast}(\Omega M)=\pi_{\ast}^{S}((\Omega
M)^{+}) \to {\bf Z}$.  It is clear that by changing the coefficients
of the complex $C^{\bf N}_{N}(f)$ by means of $r$ one obtains the classical Morse
complex of $f$.
\end{remark}

Assume $(f, N)$ is a local Morse function and fix an index pair $(N_{1},N_{0})$ for
$I_{N}(f)$. As in classical Morse theory \cite{Mi},  $N_{1}$ is obtained from $N_{0}$ by the
succesive attachement of precisely one $k$-handle for each critical point of
$f|_{N}$ of index $k$. Of course, the Conley index of $I_{N}(f)$,
$c_{\gamma(f)}(I_{N}(f))$,  which will be further denoted by $c_{N}(f)$,  also admits a cell
decomposition with one cell of dimension
$k$ for each point in $Crit^{N}_{k}(f)$. Moreover, the integral homology of the standard
Morse complex of $f$ verifies $H_{\ast}(C^{Mo}(f))\approx
\overline{H}_{\ast}(c_{N}(f);{\bf Z})$. For further use, fix also the notation
$\overline{c}_{N}(f)$ for the global Conley index of $I_{N}(f)$.

If $P,Q$ are consecutive critical points of $f$ in $N$ denote by $I(P,Q)$
the isolated invariant set consisting of the union of the points $P$, $Q$ and all the points
situated on flow lines joining $P$ to $Q$ and contained in $N$.
The triple $(I(P,Q);P,Q)$ is an attractor
repellor pair.
There is a triple $(N''',N'',N')$ of compact  subsets in $N$ such that
$(N''',N')$ is an index pair for $I(P,Q)$, $(N''',N'')$ is an index pair of $P$
and $(N'',N')$ is an index pair for $Q$. Moreover, if $p=ind(P)$,
$q=ind(Q)$ there are cofibration sequences $S^{q-1}\to N'\to N''$
and $S^{p-1}\to N''\to N'''$. By extending the first cofibration sequence one step
to the right and composing with the attaching map of $S^{p-1}$  one obtains
a map $\delta: S^{p-1}\to N'\to S^{q}$ called the relative attaching map of $P$ and $Q$.
It is easy to see that we have $\delta(P,Q)=\Sigma \delta$ where $\delta(P,Q)$ is the
connection map coming out of the attractor-repellor cofibration sequence associated to
$(I(P,Q);P,Q)$. It was shown by Franks \cite{Fr} that
$\delta$ coincides up to sign with the Thom map constructed from the
framed embedding $Z(P,Q)\subset W^{u}_{N}(P)\bigcap f^{-1}(a)\approx S^{p-1}$.

\subsection{Cofibrations and relative Hopf invariants.} \label{subsec:cof}

\

This sub-section is a digression into some elementary homotopy theory.
We will work here in the pointed category of spaces with the homotopy type of
$CW$-complexes. We will assume some familiarity with the language of (homotopy)
fibrations, cofibrations, push-outs and pull-backs as found for example in \cite{Ma}.

In particular, a sequence of two maps $A\stackrel{i}{\to} X\stackrel{i'}{\to} Y$ is a
cofibration sequence if $Y=CA\cup_{i}X$ where $CA$ is the reduced cone on $A$ and the
second map is the inclusion $X\hookrightarrow CA\cup_{i}X$.
Given such a cofibration
sequence,  there is a co-action map $\nabla:Y\to  \Sigma A\vee Y$ defined as the projection of
$Y=CA\bigcup_{A} X$  onto $Y/(A\times \{1/2\})$  where $A\times
\{1/2\}\subset CA=(A\times [0,1])/(A\times \{1\}\bigcup \ast\times [0,1])$. Assume that there
is a second cofibration sequence $\Sigma B\stackrel{j}{\to} Y\stackrel{v}{\to} Z$.
We first recall that the composition
$\Sigma B\stackrel{j}{\to}Y\stackrel{\nabla}{\to}\Sigma A\vee Y\to\Sigma A$
where the last map is the projection on the first term is called the relative attaching map
of $\Sigma B$ and $A$ (it generalizes the map $\delta$ introduced at the end of the last paragraph).

Consider the projection on the second term $p:\Sigma A\vee Z  \to Z$. The homotopy
fibre of this map is homotopy equivalent to $\Sigma A\wedge (\Omega Z ^{+})$ where
$-^{+}$ indicates the disjoint union with a point (the precise homotopy equivalence that is
used follows from diagram (\ref{equ:cube2}) below). It is easy to see that $\Omega k$  admits
a homotopical retraction. This implies, by adjunction, that the inclusion
of this fibre $k: \Sigma A\wedge (\Omega Z ^{+})\to \Sigma A\vee Z$ has the following
property:  if two maps out of a suspension and with values in the domain of $k$ are
homotopic after composition with $k$, then the two maps are homotopic.

Let $u$ be the composition
$u:\Sigma B\stackrel{j}{\to} Y\stackrel{\nabla}{\to}\Sigma A\vee
 Y\stackrel{ id\vee v}{\to}\Sigma A\vee Z$ . Then $p\circ u$ is null-homotopic and,
in view of what was said above, it has a lift $u':\Sigma B\to \Sigma A\wedge (\Omega Z
^{+})$ that is unique up to homotopy.

\begin{definition}\label{hopf} In the setting above, the Hopf invariant of
$j$ relative to $i$, $H(j,i)$, is the homotopy class of $u'$ in $[\Sigma B, \Sigma A\wedge
(\Omega Z^{+})]$
\end{definition}

This version of relative Hopf invariants is slightly different from that used in
\cite{Co2}. The relation between the two is that the projection of $H(j,i)$
on the factor $\Sigma A\wedge \Omega Z$ of $\Sigma A\wedge \Omega
X^{+}=\Sigma A\vee\Sigma A\wedge \Omega Z$ equals the Hopf invariant from
\cite{Co2}.

I am indebted to Bill Richter who has first suggested that the definition above can be
of use in this paper.

\section{Hopf invariants and connecting manifolds.}\label{sec:conn}

\

The simplest type of attractor-repellor pair is that formed by two non-degenerate consecutive
critical points $P,Q$ in a Morse-Smale gradient  flow together with the flow lines joining them.
In this case the attractor and repellor are completely understood  and the only
question is to "estimate" the connecting manifold $Z(P,Q)$. Part of the answer is contained
in the result of Franks \cite{Fr} which provides an estimate of the framed cobordism Thom
map of $Z(P,Q)$ (see above  \S \ref{subsec:rec}). The presence of Hopf
invariants in this context is natural if one views them as a homotopical way to
measure the ``strength" of the binding of two cells (or cones) that are attached successively.

The main result of this section, extends the result of Franks by giving a homotopical
description of the Thom {\em bordism} map provided by the framing $Z(P,Q)\hookrightarrow
S^{u}(P)$ and  the map $j(P,Q):Z(P,Q)\to \Omega M$.  More precisely,
consider a fixed local Morse function $(f,N)$  (see Definition \ref{def:loc_m})
Assume that $P$ and $Q$ are consecutive critical points of $f$ in $N$ of indexes respectively
$p$ and $q>0$ and recall that $Z(P,Q)$ is the connecting manifold of $P$ and $Q$. As
explained in \S \ref{subsec:rec},  $Z(P,Q)$ is a manifold of dimension $p-q-1$ which is
embedded with a fixed framing in $S^{u}(P)\approx S^{p-1}$; there is also the map
$j(P,Q):Z(P,Q)\to\Omega M$. The Thom-Pontryaguin construction applied to this data
produces a homotopy class $h(P,Q)\in \pi_{p-1}(\Sigma^{q}(\Omega X^{+}))$ such that the
stable  class of $h(P,Q)=[Z(P,Q)]$.  Let $(N_{1},N_{0})$ be an index
pair of the maximal invariant set of the negative gradient flow of $f$ inside $N$, $I_{N}(f)$.
We have  $N_{1}\subset N$ and $N_{1}$ is obtained from $N_{0}$ by
succesive cell attachments corresponding to the points in $Crit^{N}(f)$.
 In particular, because $P$ and $Q$ are consecutive we have the
cofibration sequences coming from the successive passage through these two critical points
$S^{q-1}\stackrel{j_{Q}}{\to} N'\to N''$, $S^{p-1}\stackrel{j_{P}}{\to} N''\to N'''\subset
N_{1}$ (where $N'$ can be viewed as the union of $N_{0}$ and all the handles associated to
 critical points ($\in N$) appearing before $Q$ in the negative gradient flow). By a
slight abuse of terminology, we denote by $H(j_{P},j_{Q})$  the relative Hopf invariant of
$j_{P}$ and $j_{Q}$  composed with the inclusion
$\Sigma^{q}\wedge (\Omega N''')^{+}\hookrightarrow\Sigma^{q}\wedge \Omega M^{+}$.

Here is the main result of this section.
\begin{theorem}\label{theo:hopf_conn}
We have $H(j_{P},j_{Q})=h(P,Q)$.
\end{theorem}

\begin{remark}\label{rem:hopf_conn} a. We have not been explicit enough in the orientation
choices for the Hopf invariants. Therefore, the formula in the statement above
should be interpreted as equality up to sign.

b. It has been shown in \cite{Co2} that in the {\em compact context}
we have the equality in Theorem \ref{theo:hopf_conn} {\em after one suspension}.
The proof of Theorem \ref{theo:hopf_conn} below is new, more general and leads to a stronger
result. It is also more satisfying conceptually as it distinguishes better the purely homotopical
part (Proposition \ref{prop:hopf}) from the rest of the argument. Moreover, it is also
more similar in spirit to the proof of the result in \cite{Fr} which it
extends.

c.  To be able to use Theorem \ref{theo:hopf_conn} to compute homotopically the
extended Morse complex $C^{I}_{N}(f)$ one needs to know the relative cell decomposition
induced by $f$ on an index pair $(N_{1},N_{0})$ of $I_{N}(f)$. This is somewhat
unsatisfactory because considerable information on the function $f$ is already necessary to
know this decomposition. However, we will see later that, under certain minimality
assumptions, this complex can be determined  from only the homotopy type of such a pair
$(N_{1},N_{0})$.

d. There is a significant advantage in knowing the Thom bordism  map
$h(P,Q)$ with respect to only the cobordism one:  this map
allows one to relate properties of $Z(P,Q)$ to global properties of $M$.  Notice also that
this bordism map contains finer  homotopical data than that provided by Conley's connection
map associated to the attractor-repellor pair in question.
\end{remark}

The proof of Theorem \ref{theo:hopf_conn} is contained in the first two parts of this section.
The first step is purely homotopical and is contained  in \S \ref{subsec:hopf}. The
 rest of the proof is in \S \ref{subsec:geom_hopf}. The last part of the section contains
a number of applications indicating certain rigidity properties of the extended
Morse complexes.

\subsection{Hopf invariants as relative attaching maps.} \label{subsec:hopf}

\

This sub-section is of a purely homotopical nature and can be read independently
of the rest of the paper. The first paragraph contains the main homotopical
tool, Proposition \ref{prop:hopf}, needed to prove Theorem \ref{theo:hopf_conn}.
It shows that relative Hopf invariants can be read out of certain natural relative-attaching
maps. The second paragraph discusses some direct consequences of this result. In particular
we introduce Hopf complexes which are a homotopical analogue of the extended Morse
complexes. In the third paragraph we present some immediate consequences concerning
the structure of Thom spaces and the detection of non-smoothable Poincar\'e duality spaces.

We will work here in the pointed category of spaces with the
homotopy type of finite $CW$-complexes and will assume, as in \S\ref{subsec:cof},
some familiarity with the the language and the basic results in \cite{Ma}.

\subsubsection{Cone-decompositions and Hopf invariants.}\label{subsubsec:cone}

\

Consider  a pair
$(X,X_{0})$ of
$CW$-complexes. A cone-decomposition of length
$k$ of $X$ relative to $X_{0}$ is a sequence of cofibrations sequences
\begin{equation}\label{eq:conedecomp}
A_{i}\stackrel{j_{i}}\to X_{i}\to X_{i+1}
\end{equation}
 with$0\leq i < k$ such that $X_{k}\he X$ and for each $i>0$
$A_{i}$ is a suspension. Denote by $\delta(i+1,i):A_{i+1}\to \Sigma A_{i}$
the relative attaching maps obtained by composing $j_{i+1}$
with the connectant $\delta_{i}:X_{i}\to \Sigma A_{i}$.

Consider the path-loop fibration $\mathcal{P}_{X}:\Omega X\to
PX\to X$. By pulling back this fibration over the push-out squares

$$\xy\xymatrix@-4pt{
A_{i}\ar[r] \ar[d]  & X_{i} \ar[d] \\
CA_{i} \ar[r] & X_{i+1} }
\endxy $$
we obtain push-out squares \cite{Co0}

$$\xy\xymatrix@-4pt{
A_{i}\times \Omega X \ar[r] \ar[d]  & E_{i} \ar[d] \\
CA_{i}\times\Omega X\ar[r] & E_{i+1} }
\endxy $$
where $E_{i}$ is the total space of the fibration of basis $X_{i}$ obtained by pull-back over
the inclusion $X_{i}\to X$ from $\mathcal{P}_{X}$.
 By collapsing $\Omega X$ to a point
in each corner of the previous push-out we obtain a homotopy cofibration sequence
\begin{equation}\label{equ:topconedecomp}
 A_{i}\wedge(\Omega X)^{+}\stackrel{j_{i}'}{\to} E_{i}/\Omega X\to
E_{i+1}/\Omega X
\end{equation}
We have used here the formula $(\Sigma T\times K)/K\he \Sigma T\wedge K^{+}$.
Of course, we have one such cofibration sequence for each $i$, $0\leq i<k$.
Each cofibration sequence (\ref{equ:topconedecomp}) can be extended one step to the right
thus producing a connectant $\delta_{i}':E_{i+1}/\Omega X\to \Sigma A_{i}\wedge {\Omega
X}^{+}$. Let $\delta'(i+1,i)$ be the composition $j_{i+1}'\circ \delta_{i}' $. In other words,
$\delta'(i+1,i)$ is the relative attaching map coming out of two succesive cofibration
sequences (\ref{equ:topconedecomp}).

\begin{proposition}\label{prop:hopf} The map $\delta'(i+1,i)$ and the following
composition
$$\xy\xymatrix@R-20pt{
A_{i+1}\wedge \Omega X^{+}\ar[rr]^{(s\circ H(j_{i+1},j_{i}))\wedge
id\ \ \ \ }& & \Sigma A_{i}\wedge \Omega X^{+}\wedge \Omega
X^{+}\ar[rr]^{\ \ \ \ \ \ \ \ \  \ \ \ \ \ \  id\wedge \nu  }& &\\
 \ar[r] & \Sigma A_{i}\wedge\Omega X^{+}& &  &}
\endxy $$
coincide up to homotopy.  Here $s$ is the inclusion of $X_{i+1}$ into $X$
 and $\nu$ is induced by the multiplication $\Omega X\times \Omega X\to
\Omega X$.
\end{proposition}
\begin{proof}
Consider the cube
\begin{equation}\label{equ:cube1} \begin{minipage}{5in}
\xy\xymatrix@C+15pt{
A_{i}\times\Omega X \ar[rr]^{l_{i}}\ar[dd]_{p_{1}}\ar[dr] & &E_{i}\ar@{.>}[dd]\ar[dr] &
\\ &CA_{i}\times \Omega X\ar[rr]\ar[dd] & &E_{i+1}\ar[dd]\\
A_{i}\ar@{.>}[rr]\ar[dr] &&X_{i}\ar@{.>}[dr]& \\
& CA_{i} \ar[rr] & & X_{i+1} }
\endxy \end{minipage}
\end{equation}

As discussed before this has a push-out as bottom and top squares and has pull-backs
as lateral faces. It is obtained by pull-back from the path-loop fibration $\mathcal{P}_{X}$
and by making use of the inclusion $X_{i+1}\hookrightarrow X$.  Consider now the following
diagram
\begin{equation}\label{equ:dosq1}\begin{minipage}{5in}
\xy\xymatrix@C+20pt{
A_{i+1}\times \Omega X\ar[r]^{l_{i+1}}\ar[d]_{p_{1}}&E_{i+1}\ar[d]\ar[r]^{m_{i}}
 & \Sigma A_{i}\wedge\Omega X^{+} \ar[d]\\
A_{i+1}\ar[r]_{j_{i+1}}& X_{i+1}\ar[r] &\Sigma A_{i} }
\endxy \end{minipage}\end{equation}
Here the square at the left corresponds to the square at the back  of the cube
(\ref{equ:cube1}) for the index $i+1$ and the square at the right is obtained
by taking cofibers in the square at the right of (\ref{equ:cube1}) for index $i$.
The composition $m_{i}\circ l_{i+1}$ is clearly null-homotopic when restricted
to $\ast\times\Omega X\hookrightarrow A_{i+1}\times\Omega X$. Therefore,
the map $\delta'(i+1,i)$ is obtained from this composition by collapsing to
a point $\Omega X$ in the domain and the image of $l_{i+1}$. Let $r:A_{i+1}
\longrightarrow A_{i+1}\times \Omega X\stackrel{l_{i+1}}{\longrightarrow} E_{i+1}$
be such that the first map in the composition is the inclusion on the first component.
The map $l_{i+1}$ factors  in the following way
\begin{equation}\label{equ:map1}
l_{i+1}:A_{i+1}\times\Omega X\stackrel{r\times
id}{\longrightarrow}E_{i+1}\times\Omega X \stackrel{\nu'}{\longrightarrow} E_{i+1}
\end{equation}
where the map $\nu'$ is the holonomy of the homotopy  fibration
$E_{i+1}\to X_{i+1}\to X$  and is well defined up to homotopy.

We recall that for an arbitrary fibration $F\to E\to B$. the
holonomy is obtained by transforming $F\to E$ into a fibration $\mathcal{F}$ and taking
the pull-back of $\mathcal{F}$ over the  inclusion $F\hookrightarrow E$. The fact that the
composition $F\to E\to B$  is trivial implies that the total space of $\mathcal{F}$ is
homotopy equivalent to the product $F\times\Omega B$ and the holonomy
$\nu_{\mathcal{F}}:F\times\Omega B\to F$ is the resulting map. If $E$ is contractible, then
$\nu:\Omega B\times\Omega B\to \Omega B$ is, up to homotopy, the usual multiplication.

Here is the justification of the factorization in (\ref{equ:map1}). Up to
homotopy, the map $j_{i+1}$  can be written as the composition
$A_{i+1}\stackrel{r}{\longrightarrow}E_{i+1}\longrightarrow X_{i+1}$. This means
that $l_{i+1}$ is the composition of $\nu'$ with the top map in the pull-back diagram
$$\xy\xymatrix{
A_{i+1}\times\Omega X\ar[d]_{p_{1}}\ar[r]& E_{i+1}\times\Omega X\ar[d]^{p_{1}}\\
A_{i+1}\ar[r]_{r} & E_{i+1} }
\endxy $$
but this map is $r\times id$  (the diagram is obtained immediately by transforming
$E_{i+1}\to X_{i+1}$ into a fibration and pulling back over the composition
$A_{i+1}\stackrel{r}{\to}E_{i+1}\to X_{i+1}$).

The proof of the proposition is now reduced to showing the commutativity of the
following diagram
\begin{equation}\label{equ:key}\begin{minipage}{3in}
\xy\xymatrix@C+25pt@R+20pt{
A_{i+1}\times\Omega X\ar^{r\times id}[r]\ar_{s\circ H(j_{i+1},j_{i})\times id \ \ \ \ }[dr] &
E_{i+1}\times\Omega X\ar^{\nu'}[r]\ar^{m_{i}\wedge id}[d] & E_{i+1}\ar^{m_{i}}[d] \\
& \Sigma A_{i}\wedge \Omega X^{+}\wedge\Omega X^{+}\ar_{id\wedge\nu}[r]&\Sigma
A_{i}\wedge \Omega X^{+} }
\endxy\end{minipage}\end{equation}

The commutativity of the triangle on the left in (\ref{equ:key}) follows immediately from the
fact that the map $m_{i}$  appears in the following diagram in which all squares are
pull-backs (and by performing the ``holonomy" construction on the two fiber inclusions
represented by the two vertical arrows on the left of the diagram).
\begin{equation}\label{equ:ident}\begin{minipage}{5in}
\xy\xymatrix@C+25pt{
E_{i+1}\ar^{m_{i}}[r]\ar[d]&\Sigma A_{i}\wedge\Omega X^{+}\ar[d]\ar[r]&\ast\ar[d]\\
X_{i+1}\ar_{\nabla'}[r]& \Sigma A_{i}\vee X\ar[r]& X }
\endxy\end{minipage}\end{equation}
Here the map $\nabla'$ is the composition $X_{i+1}\stackrel{\nabla}{\longrightarrow}
\Sigma A_{i}\vee X_{i+1}\hookrightarrow \Sigma A_{i}\vee X$. In turn, to verify
(\ref{equ:ident}) is is enough to compare the cube below, which is obtained
by pull-back from the fibration $\mathcal{P}_{X}$, with the similar cube
(\ref{equ:cube1}).
\begin{equation}\label{equ:cube2} \begin{minipage}{5in}
\xy\xymatrix@C+15pt{
A_{i}\times\Omega X \ar[rr]\ar[dd]_{p_{1}}\ar[dr] & &PX \ar@{.>}[dd]\ar[dr] & \\
 &CA_{i}\times\Omega X\ar[rr]\ar[dd] & & F_{i}\ar[dd]\\
A_{i}\ar@{.>}[rr]\ar[dr] &  &X\ar@{.>}[dr]& \\
& CA_{i} \ar[rr] & & \Sigma A_{i}\vee X }
\endxy \end{minipage}
\end{equation}
The map $A_{i}\to X$ is the inclusion and is clearly null-homotopic. This
shows that the square on the bottom is indeed a (homotopy) push-out.
This cube also explicits the identification of the homotopy fibre $F_{i}$ of
$\Sigma A_{i}\vee  X\to X$ with $\Sigma A_{i}\wedge \Omega X^{+}$ which
is given by  the homotopy equivalence of the cofibres of the oblique maps in the top square.
The cube (\ref{equ:cube1}) maps into (\ref{equ:cube2}) the maps defined on each corner being
induced by the inclusions on the bottom. The squares at the left in each of these two cubes
being the same (and the respective map being the identity) we obtain the claim in
(\ref{equ:ident}).

We are now left with verifying the commutativity of the square in (\ref{equ:key}). In view
of (\ref{equ:ident}) we have the following commutative diagram in which the columns are
pull-backs of $\mathcal{P}_{X}$
\begin{equation}\label{equ:fibr}\begin{minipage}{5in}
\xy\xymatrix@C+20pt{
E_{i}\ar[d]\ar[r] & E_{i+1}\ar^{m_{i}}[r]\ar[d] & \Sigma A_{i}\wedge \Omega
X^{+}\ar[d]\\ X_{i}\ar[r] & X_{i+1}\ar[r] & \Sigma A_{i}\vee X }
\endxy\end{minipage}\end{equation}
As in the construction of the holonomy of an arbitrary principal fibration we may pull-back
all the fibrations in (\ref{equ:fibr}) over the projections thus getting a new commutative
square relating the various holonomies
\begin{equation}\label{equ:hol}\begin{minipage}{5in}
\xy\xymatrix@C+20pt{
E_{i}\times \Omega X\ar_{\nu'}[d]\ar[r] & E_{i+1}\times \Omega X
\ar_{\nu'}[d]\ar^{m_{i}\times id}[r] & \Sigma A_{i}\wedge\Omega X^{+}\times \Omega X
\ar^{\nu''}[d] \\
E_{i}\ar[r] &E_{i+1}\ar[r]^{m_{i}} &\Sigma A_{i}\wedge \Omega X^{+} }
\endxy \end{minipage}\end{equation}
The desired commutativity now follows by noticing that after collapsing to a point
$\ast\times\Omega X$ in the top row of (\ref{equ:hol}) we obtain a cofibration sequence
and moreover $\nu''$ gives precisely $id \wedge \nu$ after this operation. This last fact
is seen by constructing holonomies over the vertical maps in (\ref{equ:cube2}) to
express the holonomy $\nu''$.
\end{proof}

In the remainder of this subsection we discuss some direct, useful homotopical
consequences of this proposition. The proof of Theorem \ref{theo:hopf_conn} is continued
in \S\ref{subsec:geom_hopf}.

\subsubsection{Hopf complexes.}\label{subsubsec:cor_hopf}

\

Assume that $(X,X_{0})$ is a relative $CW$-complex with a fixed cell decomposition.
Suppose that $X\to Y$ is a fixed map and $Y$ is simply connected.
Denote the set of $i$-cells in the cell- decomposition of $(X,X_{0})$ by $Cell_{i}$.
Let $I=(a_{1}<a_{2}<\ldots )\subset {\bf N}$ be such that if $Cell_{i}\not=\emptyset$ then
$i\in I$. Define $(C^{I}(X,X_{0}),d^{I})$ by letting $C^{I}(X,X_{0})$ to be the free
$\pi_{\ast}^{S}(\Omega Y^{+})$-module generated by $Cell_{\ast}$ and
if $|e|=a_{i}$, $e\in Cell_{\ast}$, let
$d^{I}(e)=\sum_{f\in Cell_{a_{i-1}}} [H(e,f)]f$ with $[H(e,f)]$
the stable homotopy class of the relative Hopf invariant $H(j_{e},j_{f})$ of the attaching
map $j_{e}$ of the cell $e$ relative to the attaching map
$j_{f}$ of the cell $f$ viewed inside $\pi_{\ast}^{S}(\Omega Y^{+})$ by making use of the
map $X\to Y$. By taking a look to the definition of Hopf invariants in \S\ref{subsec:cof} we
see that a special convention is necessary for the coefficients $[H(e,f)]$ if  the dimension of the
cell $f$ is $0$. We distinguish two cases.
First, if $dim(e)=1$ we let $\delta(e,f)=$ the coefficient
of $f$ in the differential of $e$ in the cellular complex associated to the fixed
cell-decomposition of the pair $(X,X_{0})$. We
put $H(e,f)=\delta(e,f)\in\pi_{0}^{S}(\Omega Y^{+})={\bf Z}$.
If $dim(e)=p\geq 2$ the attaching map of the cell $e$ is of the form $S^{p-1}\to
X_{0}\coprod\ldots\coprod f$. We take $H(e,f)$ to be $0$ if the image of $S^{p-1}$
is disjoint from $f$ and $H(e,f)$ is the adjoint of the map $e/\partial e\to Y$
if the image of $S^{p-1}$ equals $f$.

\begin{corollary}\label{cor:hopf}With the notations above
$(C^{I}(X,X_{0}),d^{I})$ is a chain complex. A complex of this type will be called a Hopf
complex of $(X,X_{0})$ with coefficients in $\pi_{\ast}^{S}(\Omega Y^{+})$.
\end{corollary}

\begin{proof}
Apply Proposition \ref{prop:hopf} to express
the relevant  Hopf invariants in terms of relative attaching maps. Use the fact that the
composition of any two succesive relative attaching maps is null-homotopic to deduce relations
among the relevant Hopf invariants. After stabilization, these relations are precisely those
needed to show
$(d^{I})^{2}=0$.
\end{proof}

\begin{remark} a. The Hopf complexes as constructed above are closely related to
the twisted chain complexes of Baues \cite{Ba}. In fact, $(C^{I}(X,X_{0}),d^{I})$ can be
identified with a chain complex that is easely deduced from this
twisted complex. We will not explicit this  identification here. It is immediately implied by
the fact that the partial suspension of \cite{Ba} is also given by the following
construction. Assume that $A\stackrel{f}{\longrightarrow} B\vee Y$ is a map of pointed
$CW$-complexes which is null on $Y$. Then the partial suspension of $f$ is homotopic to the
composition $i \circ \Sigma f'$. Here  $f'$ is the unique lift (up to homotopy) of $f$ to the
homotopy fibre of $B\vee Y \longrightarrow Y$ and $i$ is the fibre inclusion
in  the homotopy fibration $\Sigma B\wedge \Omega Y^{+}\to \Sigma B\vee Y\to Y$.

b. Smaller  is the index set used in defining a Hopf complex $C^{I}(X,X_{0})$,  more
information is contained in the differential of this complex. Indeed, enlarging this index set
changes the respective complex by possibly truncating the differential. Moreover,
because $Y$ is simply-connected we notice that $C^{\bf N}(X,X_{0})$ is simply
the cellular complex associated to the fixed relative cell decomposition of $(X,X_{0})$
tensored with $id_{\pi_{\ast}^{S}(\Omega Y^{+})}$.
\end{remark}

It is natural to compare Hopf complexes.

Assume $(X',X)$ is a relative $CW$-complex. Fix two different
relative cell-decompositions for $(X',X)$ and a map $X'\to Y$ with $Y$ simply connected.
Let  $C^{I}(X',X)$ and $C^{I}_{1}(X',X)$ be the two Hopf complexes with coefficients
in $\pi_{\ast}^{S}(\Omega Y^{+})$ respectively associated to these two cell
decompositions.

\begin{lemma}\label{lem:h_comp} There is a chain equivalence of the two cellular complexes
associated to the two cell decompositions of $(X',X)$  which extends to a chain equivalence
of $C^{I}(X',X)$ and $C^{I}_{1}(X',X)$  whenever these Hopf complexes are both
defined. In particular, if both cell decompositions are minimal in the sense that  in each
dimension $k$  they have the same number of cells and this number equals $rk
(H_{k}(X',X; {\bf Z}))$, then the Hopf complexes above are isomorphic.
If such minimal cell-decompositions exist, then the corresponding Hopf complex that  appears
for $I=I_{min}=\{k : H_{k}(X',X;{\bf Z})\not=0\}$ will be denoted by
$(\mathcal{C}_{min}(X',X),d_{min})$ and is well defined up to isomorphism.
\end{lemma}

\begin{proof}

Assume $\vee S^{k}\to X^{(k)}\to X^{(k+1)}$ and $\vee S_{1}^{k}\to
X^{(k)}_{1}\to X_{1}^{(k+1)}$ are the cofibration sequences in the two cell
decompositions.  We have that
$X^{(0)}$, $X^{(0)}_{1}$ are equal to the disjoint union of $X$  and some points.  For
some $m\in {\bf N}$ we have $X^{(m)}\simeq X^{(m)}_{1}\simeq X'$. By restricting this
homotopy equivalence to each of the skeleta we get maps $\phi_{k}:X^{(k)}\to
X_{1}^{(k)}$ that commute with the obvious inclusions and restrict to the identity on $X$.
We obtain the commutative diagram below with the top and bottom rows
cofibration sequences and where the existence of the doted arrow remains to be discussed.

\begin{equation}\label{equ:comp_h}\begin{minipage}{5in}
\xy\xymatrix{
\vee S^{k}\ar[r]\ar@{.>}_{\phi''_{k+1}}[d]& X^{(k)}\ar_{\phi_{k}}[d]\ar[r]&
X^{(k+1)}\ar^{\phi_{k+1}}[d]\ar[r]&\vee S^{k+1}\ar^{\phi'_{k+1}}[d] \\
\vee S^{k}_{1}\ar[r]& X^{(k)}_{1}\ar[r]&
X^{(k+1)}_{1}\ar[r]&\vee S^{k+1}_{1} }
\endxy \end{minipage}\end{equation}

It is a standard fact that the maps $\phi'_{t}$  induce a a chain equivalence of cellular
complexes  $\phi=H_{\ast}(\phi'_{\ast})$. There exists always a map $\phi''_{k+1}$
which is a desuspension  of $\phi'_{k+1}$.
If $k\geq 2$ this desuspension is uniquely defined and we first notice that, under this same
assumption, it does make commute the left square of the diagram. This is quite immediate
if $X'$ is simply connected and therefore both $X^{(k)}$ and $X^{(k)}_{1}$ are so.
Indeed, the composition
$\pi_{k}(\vee S^{k})\approx H_{k}(\vee S^{k})\approx
H_{k+1}(X^{(k+1)},X^{k})\approx \pi_{k+1}(X^{(k+1)},X^{(k)})\to \pi_{k}(X^{(k)})$
describes precisely the attaching maps of the $k+1$ cells in the first cell-decomposition and
the analogue statement is valid, of course, for the second decomposition. But all the
maps in these compositions commute with the map induced in relative homology by the
pair $(\phi_{k+1},\phi_{k})$. If $X'$ is not simply-connected (but still $k\geq 2$)
 one needs to
first apply the above argument to the universal covering $\tilde{X'}$
relative to the cell-decompositions obtained by lifting those given for $X'$.

Our purpose now is to show that the map of complexes $\phi^{I}:C^{I}(X',X)\to
C^{I}_{1}(X',X)$ defined by  $\phi\otimes id_{\pi_{\ast}^{S}(\Omega
Y^{+})}$  is a chain map. This is immediate if the
commutativity in the  left two squares of the diagram (\ref{equ:comp_h}) would hold for all
$k$ because all the relevant constructions (co-actions, lifts, etc) would be comparable by maps
induced by the $\phi''_{k+1}$'s.  However, in general, this commutativity is not valid for
$k\leq 1$.  We first assume that $2\in I$ and let
$q$ be the first index in $I$ that is bigger than $2$. It is easy to see, due to the relation
between Hopf complexes and cellular complexes, that  for dimensions $t\leq 2$
 the maps $\phi^{I}$ commutes with the differentials
in the respective Hopf complexes. This also happens for $t>q$, as indicated above, because of
the commutativity of the diagram (\ref{equ:comp_h}). Moreover, the commutativity
in the left square of (\ref{equ:comp_h}) is valid for $k=q-1$. In constructing the
relative Hopf invariants of a cell of dimension $q$ and one of dimension $2$
, besides the attaching maps for the cell of dimension $q$, we also need
the maps $X^{(2)}\to X^{(2)}\vee (\vee S^{2})\to Y\vee (\vee S^{2})$ where the first map in
this composition is the co-action and the last is induced by inclusion. It is easy to see that the
remaining commutativity necessary for
$\phi^{I}$  to be a chain morphism reduces to showing that the
square below commutes.

\begin{equation}\label{equ:comp_sq}\begin{minipage}{5in}
\xy\xymatrix{
X^{(2})\ar_{\phi_{2}}[d]\ar[r]& Y\vee (\vee S^{2})\ar^{id\vee \phi'_{2}}[d]\\
X^{(2)}_{1}\ar[r] & Y\vee (\vee S^{2}_{1}) }
\endxy \end{minipage}\end{equation}

Because $Y$ is simply connected, this commutativity is satisfied if its homology
version is true which is immediate.
The same argument applies to the case when $2\not\in I$ but $1\in I$ and shows that
$\phi^{I}$ is a chain morphism.

Because $\phi$ is a chain equivalence $\phi^{I}$ is forced to
be one also. The reason is that the differential in the cellular complex associated
to a cell decomposition of $(X',X)$ is contained in the differential of
each Hopf complex $C^{I}(X',X)$ which is defined.

The last claim
in the statement is obvious because for a minimal cell decomposition the
cellular differential is trivial and therefore, as the comparison morphism $\phi^{I}$
is a chain equivalence, it is in fact an isomorphism.
\end{proof}

\begin{remark} \label{lem:u_hopf} A condition that insures the existence of a minimal
cell decomposition as in the lemma above is the following: $(X',X)$ has torsion free integral
homology and both $X'$ and $X$ are simply-connected.
\end{remark}

Here is a relation between relative Hopf invariants and a certain
Atiyah-Hirzeburch spectral sequence.  We will only formulate this consequence in the
setting of the Remark \ref{lem:u_hopf}.

Therefore, $(X,X_{0})$ is a $CW$-pair with both
$X_{0}$ and $X$ simply-connected and such that the integral homology of the pair is torsion
free. We consider Hopf complexes with coefficients in $\pi_{\ast}(\Omega X^{+})$ (in other words
we take $Y=X$ and the relevant map $X\to Y$ is the identity).
We have a ring map  $h:\pi_{\ast}^{S}(\Omega X^{+})\to H_{\ast}(\Omega X^{+};{\bf
Z})$. We can change the coefficents in $\mathcal{C}_{\min}(X,X_{0})$ by using this map and we
denote by $(\check{\mathcal{C}}_{min}(X,X_{0}),\check{d}_{min})$ the resulting chain
complex of $H_{\ast}(\Omega X^{+};{\bf Z})$-modules.  Let
$(\mathcal{E}^{r}_{pq},D^{r})$ be the Atiyah-Hirzebruch-Serre spectral sequence of the
relative fibration $\mathcal{P}'$ of basis
$(X,X_{0})$ induced from the fibration $\Omega X\to PX\to X$ with
 $$\mathcal{E}^{2}=H_{\ast}(X,X_{0};\pi_{\odot}^{S}(\Omega X^{+}))\approx
H_{\ast}(X,X_{0};{\bf Z})\otimes \pi_{\odot}^{S}(\Omega X^{+})$$
and which
converges to $\pi_{\ast+\odot}^{S}$. Similarly, we let $(E_{pq}^{r},d^{r})$ be the homology
Serre spectral sequence of the fibration $\mathcal{P}'$.

\begin{corollary}\label{cor:s_ss} In the setting above assume that $a,b\in{\bf N}$ are
such that $H_{k}(X,X_{0};{\bf Z})=0$ for $a<k<b$, $b-a>1$. If $H_{k}(X,X_{0};{\bf Z})\not=0$ for
$k\in\{a,b\}$, then
\begin{itemize}
\item[(i)] $d_{min}|_{Cell_{b}}\to \pi_{\ast}^{S}(\Omega X^{+})<Cell_{a}>$
is identified to $D^{a-b}|_{\mathcal{E}^{b-a}_{a0}}$.
\item[(ii)] $\check{d}_{min}|_{Cell_{b}}\to H_{\ast}(\Omega
X^{+};{\bf Z})<Cell_{a}>$ is identified to $d^{a-b}|_{E^{b-a}_{a0}}$.
\end{itemize}
\end{corollary}

\begin{proof} This follows by noting that the elements in $\mathcal{E}^{2}_{a0}$
(repspectively in $E^{2}_{a0}$ survive to $\mathcal{E}^{a-b}_{a0}$ (resp.
$E^{a-b}_{a0}$) for dimensionality reasons and then by recalling the definition
of the two spectral sequences considered (as described for example in \cite{Sw})
and applying the proposition.
\end{proof}

\begin{remark} The second point of the corollary is particularly useful because
it allows one to compute the homology version of the minimal Hopf complex only in terms
of the Serre spectral sequence of the path-loop fibration of $X$.  In turn, this
spectral sequence is generally very well understood. If we are interested in a non minimal
Hopf complex or even one Hopf complex that comes from a non-minimal cell decomposition
 it is still possible to estimate it, in a rather obvious way, by using these spectral sequences but
we will  not pursue further this issue here.
\end{remark}

\subsubsection{Duality.}

\

Assume that $\xi:X\longrightarrow BG$ is the classifying map of a fiber bundle of rank
$m$, of basis $X$ and of structural group $G$. Let $E^{\xi}_{S}(X)$ be the total space of the
associated spherical fibration and let $T^{\xi}(X)$ be the respective Thom space. Fix also the
cone-decomposition  $A_{i}\stackrel{j_{i}}{\to} X_{i}\to X_{i+1}$, $0\leq i<k$ as in
(\ref{eq:conedecomp}) together with the relative attaching maps $\delta(i+1,i):A_{i+1}\to
\Sigma A_{i}$.

\begin{corollary}\label{cor:cdec_thom} The Thom space $T^{\xi}(X)$ admits
a cone decomposition
\begin{equation}\label{eq:cdec_thom}
\Sigma^{m}A_{i}\stackrel{j_{i}''}{\to}T^{\xi}(X_{i})\to T^{\xi}(X_{i+1})
\end{equation}
such that the corresponding attaching maps $\delta''(i+1,i):\Sigma^{m}A_{i+1}\to
\Sigma^{m+1}A_{i}$  verify  (up to sign) $\delta''(i+1,i)=\Sigma^{m}\delta(i+1,i)+
\overline{\nu}\circ\Omega\xi\circ\Sigma^{m}H(j_{i+1},j_{i})$.
Here $\overline{\nu}:\Sigma A_{i}\wedge S^{m-1}\wedge G\to \Sigma
A_{i}\wedge S^{m-1}$ is induced by the action of $G$ on $S^{m-1}$ and $\Omega\xi
: \Sigma A_{i}\wedge S^{m-1}\wedge\Omega X\to \Sigma A_{i}\wedge S^{m-1}\wedge G$.
\end{corollary}

\begin{proof}
Notice that the exact same constructions as in the proof of the proposition can be
performed for the fibration $\xi$ instead of the fibration $\mathcal{P}_{X}$.
Of course, as $\mathcal{P}_{X}$ is universal, the map $ \Omega X\to S^{m-1}$
induces comparison maps
between the diagrams obtained for $\mathcal{P}_{X}$ and those for $\xi$.
Therefore, the argument in the proposition
describes also the relative attaching maps coming from a cone-decomposition
$A_{i}\wedge (S^{m-1})^{+}\to E^{\xi}_{S}(X_{i})/S^{m-1}\to
E^{\xi}_{S}(X_{i})/S^{m-1}$. It is now immediate to see that, by collapsing each such
cofibration sequence into the cofibration $A_{i}\to X_{i}\to X_{i+1}$ one gets
cofibration sequences $A_{i}\wedge S^{m}\to T^{\xi}(X_{i})\vee S^{m}\to
T^{\xi}(X_{i+1})\vee S^{m}$. We then collapse to a point the sphere $S^{m}$ in the
last two terms and obtain the cofibration sequences of the statement
and the claimed formula for the relative attaching maps.
\end{proof}

\begin{remark}\label{rem:j_hom} a. The formula in this corollary also appears in Dula in
\cite{Du}.

b. Of particular interest is the case when the spaces $A_{i}$ are spheres
of dimensions, respectively, $a_{i}$.
Then $H(j_{i+1},j_{i})\in \pi_{a_{i+1}}(\Sigma^{a_{i}+1}(\Omega X)^{+})$
and the stable difference $\delta''(i+1,i)-\delta(i+1,i)$  becomes equal to
$J_{G}^{a_{i}}(\Omega\xi \circ H(j_{i+1},j_{i}))$ where
$J_{G}^{a_{i}}:\pi_{a_{i+1}-a_{i}-1}\Omega^{a_{i}}\Sigma^{a_{i}}
G\to\pi^{S}_{a_{i+1}-a_{i}-1}$ is the factor of the $J_{G}$ homomorphism (as
defined for example in \cite{Wh} and also see below) $J_{G}:\pi_{\ast}G\to \pi^{S}_{\ast}$.

c. The statement of Corollary \ref{cor:cdec_thom} is clearly also valid if applied
directly to  spherical fibrations (and not only to bundles).
\end{remark}

Assume now that $X$ is a Poincar\'e duality space and let $\zeta$  be
its Spivak stable normal bundle \cite{Sp}.

We have a ring map $s:\pi_{\ast}^{S}(\Omega X^{+})\to
\pi_{\ast}^{S}$. Assume also that $X$ is simply-connected and does not have any torsion in
homology. Then $(X,\ast)$ has a minimal cell-decomposition in the sense of
the lemma above. Denote by $\mathcal{C}_{min}(X)$  the Hopf complex associated to
such a minimal cell-decomposition with coefficents in $\pi_{\ast}^{S}(\Omega X^{+})$.
Let $(\overline{\mathcal{C}}_{min}(X),\overline{d})$ be the chain complex with coefficients
in $\pi_{\ast}^{S}$ obtained from
$\mathcal{C}_{min}(X)$ by changing the coefficients via $s$. The coefficients of the
differential in this chain complex are homotopy classes of relative attaching maps of
consecutive pairs of cells.
Finally, let
$(\overline{\mathcal{C}}'_{min}(X), \overline{d}')$ be the following $\pi_{\ast}^{S}$-chain
complex. As $\pi_{\ast}^{S}$-modules we have isomorphisms
$(\overline{\mathcal{C}}'_{min}(X))_{n-k}\approx
(\overline{\mathcal{C}}_{min}(X))_{k}$  and if $x$ is a generator of
$(\overline{\mathcal{C}}_{min}(X))_{k}$ we denote by $x^{\ast}$ the corresponding
generator of $(\overline{\mathcal{C}}'_{min}(X))_{n-k}$; by definition, the differential verifies
$\overline{d}'(x^{\ast})=\sum \overline{d}'(x^{\ast},y^{\ast})y^{\ast}$
with $\overline{d}'(x^{\ast},y^{\ast})=\overline{d}(y,x)+
J_{BF}^{|x|}(\Omega\zeta \circ H(y,x))$. Here $y$ and $x$ represent cells in the minimal
cell-decomposition that have consecutive indexes in $\{k\in {\bf N}: H_{k}(X:{\bf
Z})\not=0\}$; $H(y,x)$ is the respective relative Hopf invariant;
$BF$ is the classifying space of spherical fibrations (see also Remark \ref{rem:j_hom}) and
$\overline{d}(y)=\sum \overline{d}(y,x)x$. We also recall the definition of the
$J$-homomorpphism \cite{Wh}. This is a homomorphism
$J_{SO}:\pi_{k}({\bf SO})\to\pi_{k}^{S}$ that associates to the homotopy class of a map $a
:S^{k}\to SO$ the homotopy class of a map obtained as follows. First, for $q$ sufficiently big,
let $b:S^{q-1}\times S^{k}\stackrel{id\times a}{\to} S^{q-1}\times {\bf
SO}(q)\stackrel{\mu}{\to} S^{q-1}$ where $\mu$ is the action. Define $J_{SO}(a)$ to be the
stable homotopy class of the restriction of $\Sigma b$ to
$S^{q+k}\subset \Sigma (S^{q-1}\times S^{k})$.  The morphism $J_{SO}$ has
a stable target. Because of that it factors as $\pi_{k}({\bf SO})\to \pi_{k}(\Omega \Sigma
{\bf SO})\ldots \pi_{k}(\Omega ^{q}\Sigma^{q}{\bf SO})\to \ldots \pi_{k}^{S}({\bf
SO})\stackrel{J'_{SO}}{\to}\pi_{k}^{S}$. We denote by $J^{q}_{SO}$ the
restriction of $J'_{SO}$  to $\Omega ^{q}\Sigma^{q}{\bf SO}$. The same morphism
can be defined in an analogue fashion for other groups $G$ that have the property to
act on homotopy spheres. The corresponding homomorphisms are denoted by
$J_{G}$ and $J^{q}_{G}$.

\begin{corollary}\label{cor:PD} In the setting above ($X$ simply-connected and whithout
torsion in homology) the two $\pi_{\ast}^{S}$ -chain complexes
$(\overline{\mathcal{C}}'_{min}(X), \overline{d}')$ and
$(\overline{\mathcal{C}}_{min}(X), \overline{d})$ are isomorphic. If the stable spherical
normal fibration of $X$ admits a reduction to $BG$ (for example, $BPL$ or $BSO$ etc.),
then
$$\overline{d}'(x^{\ast},y^{\ast})=+/-\overline{d}(y,x)\ mod \ Im(J_{BG}^{|x|})$$
\end{corollary}

\begin{proof}
This follows from Corollary \ref{cor:cdec_thom} and Remark \ref{rem:j_hom}.
Indeed, the Thom space of the Spivak normal stable bundle of $X$ is Spanier-Whitehead dual
to $X$ itself. Therefore, the fixed  minimal cell decomposition on $X$
induces by duality a minimal decomposition for this Thom space.
As Spanier-Whitehead duality for maps between spheres
coincides with stable equality up to sign the relative attaching maps of this decomposition
coincide stably with those of the one fixed on $X$. Moreover, as both this decomposition
and that provided by Corollary \ref{cor:cdec_thom} are minimal we obtain that the
two complexes in question are isomorphic.  The second part of the statement is clear
given the definition of the $J$-homomorphism.
\end{proof}

\begin{remark} The problem of deciding whether a given Poincar\'e duality space $X$
admits a PL or smooth structure is one of the central problems of  differential topology.
Surgery theory provides a solution to this question by reducing it to a
two stage process: first, one has to decide whether the Spivak stable normal bundle
admits a lift to $BPL$ or, respectively, $BSO$ and, if this is the case, a secondary
numerical invariant needs to vanish. The existence of a lift is theoretically the simple
step as it reduces to a typical obstruction theory problem. However, in practice,
even for spaces $X$ having a small number of cells it is not easy to show that
such a lift does not exist. Indeed, the canonical way to produce a lift is to extend
it inductively over each skeleton of $X$ and if some non-vanishing obstruction appears
at some step in the process one is forced to start back from the bottom. This direct approach has
the disadvantage that, to conclude the non-existence of a lift, one has to show that each
possible attempt leads to a non-vanishing obstruction. The corollary above provides a shortcut
to this process: if for some $x,y$ as in the corollary we do not have the relation
$\overline{d}'(x^{\ast},y^{\ast})=+/-\overline{d}(y,x)\ mod \ Im(J_{SO}^{|x|})$ for example,
then $X$ is not smoothable. The first such example of a Poincar\'e duality space $X$  with
$|x|=0$ was described by Smith in \cite{Sm} (his proof being different from the
method indicated here).  Examples with $|x|>0$ have been first described in \cite{Co2}
the proof given there was however not completely homotopical.
To the author's knowledge, there is no other general method to construct explicitely
non-smoothable Poincar\'e duality, $1$-connected, $CW$-complexes that remain
non-smoothable after  product with an arbitrary sphere.
\end{remark}

\subsection{Relating Hopf invariants and bordism maps.} \label{subsec:geom_hopf}

\

The purpose of this sub-section is to relate relative Hopf invariants to the Thom bordism
maps associated to connecting manifolds and in this way prove Theorem
\ref{theo:hopf_conn}.

We now return to the setting of this theorem.

As in \S\ref{subsubsec:cone}  let $\mathcal{P}_{M}$ be the
path-loop fibration $\Omega M\to PM\to M$. Of course, we have inclusions
$N_{0}\hookrightarrow N'\hookrightarrow N''\hookrightarrow N'''\hookrightarrow
N_{1}\hookrightarrow M$. For any subset $A\subset M$
we can pull-back $\mathcal{P}_{M}$ over $A$ thus getting a new fibration of
basis $A$ whose total space will be denoted by $T(A)$.
As in
(\ref{equ:topconedecomp}) we  obtain two cofibration sequences
$S^{q-1}\wedge \Omega M^{+}\to T(N')/\Omega M\to T(N'')/\Omega M$
and $S^{p-1}\wedge\Omega M^{+}\to T(N'')/\Omega M\to T(N''')/\Omega M$
that produce a relative attaching map $\delta:S^{p-1}\wedge\Omega M^{+}\to
S^{q}\wedge \Omega M^{+}$.

In view of Proposition \ref{prop:hopf} the statement of the theorem is clear if
we show that $\delta$ is homotopic to the composition
\begin{equation}\label{equ:TP_comp}
\xy\xymatrix@C+20pt{
S^{p-1}\wedge\Omega M^{+}\ar^{h(P,Q)\wedge
id \ \ \ \ \ \ }[r] &S^{q}\wedge\Omega M^{+}\wedge \Omega
M^{+}\ar^{\ \ \ \ \ \ id\wedge\nu}[r] & S^{q}\wedge \Omega M^{+} }
\endxy\end{equation}
where, as in Subsection \ref{subsec:hopf}, $\nu$ is the multiplication on $\Omega M$.

We denote by $\gamma$ the flow $\gamma(f)$.
The first step of the proof is to make some explicit constructions which will allow
us to express $\delta$ in terms of some simple index pair of the isolated invariant set $I(P,Q)$
which, we recall, is the union of $P$, $Q$ and of all the points situated on flow lines of
$\gamma$ that join $P$ to $Q$ and are included in $N$.

\subsubsection{A simple index pair of $I(P,Q)$.}

\

Let
\begin{eqnarray}\label{equ:cyl}
U(P;\epsilon,\tau)=\{ x\in M : \exists\ \  T_{x}\in {\bf R}\cup\{-\infty, \infty \},
\gamma_{T_{x}}(x)\in  \\
\in D(P,\tau)\cap f^{-1}(f(P)), \ \ |f(x)-f(P)|\leq\epsilon \} \nonumber
\end{eqnarray}
where $D(P,\tau)$ is the closed disk of radius $\tau$ and center $P$.
Clearly, for $\epsilon$ and $\tau$ small the set $U(P;\epsilon,\tau)$
is a manifold (with a boundary with corners) homeomorphic
to the $n$-dimensional closed disk. Let $A(P;\epsilon,\tau)=U(P;\epsilon,\tau)\cap
f^{-1}(f(P)-\epsilon)$, $B(P;\epsilon,\tau)=U(P;\epsilon,\tau)\cap f^{-1}(f(P)+\epsilon)$.
Then $A(P;\epsilon,\tau)$ is diffeomorphic to $S^{p-1}\times D^{n-p}$ and
$B(P;\epsilon,\tau)$ is diffeomorphic to $D^{p}\times S^{n-p-1}$.
In these identifications $S^{p-1}\times\{0\}$ and $\{0\}\times S^{n-p}$ correspond
respectively to the unstable and the stable manifolds of $P$ intersected respectively with
$f^{-1}(f(P)-\epsilon)$ and $f^{-1}(f(P)+\epsilon)$ .
We also denote by $A'(P;\epsilon,\tau)$  the closure of the set $\partial
U(P;\epsilon,\tau)-B(P;\epsilon,\tau)$. The two diffeomorphisms above can be made quite
explicit once Morse coordinates are chosen around $P$. Using such coordinates $(x,y)$ with
$f$ given by $-||x||^{2}+||y||^{2}+f(P)$  there is a diffeomorphism of
$A(P;\epsilon,\tau)$  with a set of the type $\{ (x,y)\in{\bf R}^{p}\times {\bf R}^{n-p}:
\epsilon\leq ||x||^{2}\leq \epsilon + \eta \ , \  -||x||^{2}+||y||^{2}=-\epsilon \}$. and a similar
identification is valid for $B(P;\epsilon,\tau)$. Moreover, the  flow induces
a diffeomorphism between $B(P;\epsilon,\tau)-W^{s}_{N}(P)$ and
$A(P;\epsilon,\tau)-W^{u}_{N}(P)$ which corresponds to the obvious twist map
$(S^{p-1}\times (0,1])\times S^{n-p-1}\longrightarrow S^{p-1}\times (S^{n-p-1}\times
(0,1])$. We also need to notice that, once a Morse chart as above is fixed, there are $n$-linearly
independent sections of the tangent space of $U(P;\epsilon,\tau)$ each in the direction
of one of the coordinates $x_{i}$ or $y_{j}$. These also provide normal framings
to the sets $Int(W^{u}_{N}(P))$ and $Int(W^{s}_{N}(P))$. These framings have the property
to  be coherent with the decomposition $A(P;\epsilon,\tau)\approx S^{p-1}\times D^{n-p}$
(and similarly for $B(P;\epsilon,\tau)$ ) in the sense that they induce the standard framing of
$S^{p-1}\times\{0\}\subset S^{p-1}\times D^{n-p}$.
In particular, we may see $A(P;\epsilon,\tau)$ as the total space of a tubular neighbourhood
of $S^{u}(P)$. We denote by $\iota_{P}^{-}$ and $\iota_{P}^{+}$ the framings of
$S^{u}(P)$ in $A(P;\epsilon,\tau)$ and respectively that of $S^{s}(P)$ in
$B(P;\epsilon,\tau)$. Of course, these framings depend on the choice of a Morse chart
but we suppose from now on that such a choice has been made for all critical points
of $f$. A different choice will not modify the isomorphism class of these framings
(except possibly by a sign).

An analogue construction is of course possible for the critical point $Q$. As the critical points
$P$ and $Q$ are consecutive we may assume (after possibly slightly modifying $f$)
that for some choice of $\epsilon$ and $\tau$  the sets $U(Q)=U(Q;\epsilon,\tau)$,
$U(P)=U(P;\epsilon,\tau)$  and $A(P)=A(P;\epsilon,\tau)$ verify
$U(P)\cap U(Q)=\{x\in A(P):\ \gamma_{\bf R}(x)\cap U(Q)\not=\emptyset\}$
(this implies, in particular, that $\epsilon$ is such that $f(Q)+\epsilon=f(P)-\epsilon$).
Moreover, the intersection of $S^{p-1}\times\{0\}\subset A(P)$ and $\{0\}\times
S^{n-q-1}\subset B(Q)$  are transverse in $f^{-1}(\epsilon+f(Q))$ (where $B(Q)=
B(Q;\epsilon,\tau)$ and we will use similar notations for the rest of the sets $A,B,A'$
associated to $U(P)$ and $U(Q)$). This means that by assuming $\tau$ sufficiently
small, the intersection of $A(P)$ and $B(Q)$ is also transverse and, by using
the identifications described before for these two sets and letting $K=U(P)\cap U(Q)$,
we have $K=A(P)\cap B(Q)\approx D^{q}\times Z(P,Q)\times D^{n-p}$.  Clearly,
$D^{q}\times Z(P,Q)\times
\{0\}\subset S^{u}(P)$ and $Z(P,Q)\hookrightarrow S^{u}(P)$ is normally framed via the
framing obtained by projection from $\iota_{Q}^{+}$.  We will denote this framing
of $Z(P,Q)$ by $\iota^{+}$ and we let $W(Z)=D^{q}\times Z(P,Q)\times
\{0\}\subset S^{u}(P)$ be the tubular neighbourhood of $Z(P,Q)$  described above.
Of course, $\iota^{+}$ identifies $W(Z)$ with the total space of the normal bundle
of $Z(P,Q)$ in $S^{u}(P)$. The framing $\iota^{+}$ is precisely the one used to
construct the bordism class of $Z(P,Q)$.

Let $W'=A'(Q)\cup (A(P)-K)$, $W''=U(Q)\cup A(P)$, $W'''=U(Q)\cup U(P)$.
These three sets are related by cofibration sequences described before the statement
of the theorem. They also verify the needed index pairs conditions and it is immediate to see
that they can be used in the constructions made at the beginning of the proof by
replacing $N',N'',N'''$ by $W',W''$ and
$W'''$, respectively. Moreover, $\delta$ is also given as the relative attaching map described in
that construction now applied to the spaces $T(W')/\Omega(M)$, $T(W'')/\Omega M$ and
$T(W''')/\Omega M$. Indeed, $(W''', W')$  is an index pair of
$I(P,Q)$ and for any other index pair $(N''',N')$ of this invariant set
we may choose $s$, $\tau$ sufficiently small such that the set $(U(Q;\epsilon,\tau)\cup
U(P;\epsilon,\tau))\bigcap f^{-1}([f(Q)-s,f(P)+s])$  (which is diffeomorphic to $U(P)\cup
U(Q)$ and has analogue properties) is contained in $N'''$. This allows one to identify the
attaching  map constructed in $T(N''')$ with that constructed in $T(W''')$.

\subsubsection{The Thom-Pontryaguin construction and $\delta$.}

\

Our next purpose  is to sufficiently explicit $\delta$ as to show that it is given by
the composition in (\ref{equ:TP_comp}). Recall that for any subspace $A\hookrightarrow
M$ we denote by $T(M)$ the total space of the fibration of basis $A$ obtained by pull-back
from $\mathcal{P}_{M}$.  By revisiting (\ref{equ:cube1}) (with $A_{i+1}=S^{u}(P)$
and $A_{i}=S^{u}(Q)$) we see that
$\delta$ is determined by the following composition
\begin{equation}\label{equ:tot_at}
\delta':T(S^{u}(P))\hookrightarrow T(A(P))\hookrightarrow
T(W'')\stackrel{l}{\longrightarrow} T(W'')/T(W')
\end{equation}
where the first three maps are inclusions and $l:T(W'')\to T(W'')/T(W')$ is the projection.
Indeed, $\delta'$ is immediately seen to be equal to the composition
$T(S^{u}(P))\to T(S^{u}(P))/\Omega M\stackrel{\delta}{\to}(T(W'')/\Omega
M)/(T(W')/\Omega M)= T(W'')/T(W')$.
At the same time $T(W'')/T(W')=T(U(Q))/T(A'(Q))$. It follows that $\delta'$ is also
equal to the composition
\begin{equation}\label{equ:tot_th}
T(S^{u}(P))\stackrel{k}{\to} T(W(Z))/T(\partial W(Z))\hookrightarrow T(U(Q))/T(A'(Q)
\end{equation}
where the map $k$ is projection on the respective quotient and the last map
is well defined because $\partial W(Z)\subset A'(Q)$.
The fibrations producing the total spaces in (\ref{equ:tot_th}) are all trivial because,
one one hand, $A'(Q)\subset U(Q)\he\ast$ and on the other $\partial W(Z)\subset W(Z)\subset
S^{u}(P)\subset W^{u}(P)\he\ast$. To get to the expression in (\ref{equ:TP_comp})
we need to explicitely use the respective trivializations.   Fix a trivialization of the fibration
$\mathcal{F}_{W^{u}(P)}$ obtained as a pull-back of $\mathcal{P}_{M}$. It produces
the following diagram where all horizontal arrows are the obvious inclusions and
the vertical arrows are homeomorphisms
\begin{equation}\label{equ:1_triv}\begin{minipage}{5in}
\xy\xymatrix@C-10pt{
\partial W(Z)\times\Omega M\ar[r]\ar[d]_{a} &W(Z)\times\Omega M\ar[d]_{a}\ar[r]
&S^{u}(P)\times\Omega M\ar[d]\ar[r] &W^{u}(P)\times \Omega M\ar[d] \\
T(\partial W(Z))\ar[r] &T(W(Z))\ar[r] & T(S^{u}(P))\ar[r] &T(W^{u}(P)) }
\endxy\end{minipage}\end{equation}

Notice that $W(Z)\subset B(Q)\subset U(Q)$. Therefore, if we fix a trivialization for the
fibration induced over $U(Q)$ (from $\mathcal{P}_{M}$) we get another similar diagram
\begin{equation}\label{equ:2_triv}\begin{minipage}{5in}
\xy\xymatrix@C-10pt{
\partial W(Z)\times\Omega M\ar[r]\ar[d]_{b} &W(Z)\times\Omega M\ar[d]_{b}\ar[r]
&B(Q)\times\Omega M\ar[d]\ar[r] &U(Q)\times \Omega M\ar[d] \\
T(\partial W(Z))\ar[r] &T(W(Z))\ar[r] & T(B(Q))\ar[r] &T(U(Q)) }
\endxy\end{minipage}\end{equation}

We may consider the composition $c=a\circ b^{-1}$ of the two homeomorphisms $a$ and
$b$ appearing  respectively in the diagrams (\ref{equ:1_triv}) and (\ref{equ:2_triv}).
It makes commutative the diagram
$$\xy\xymatrix{
W(Z)\times\Omega M\ar^{c}[r]\ar_{p_{1}}[d] & W(Z)\times\Omega M \ar^{p_{1}}[d]\\
W(Z)\ar_{id}[r] & W(Z) }
\endxy $$
however the key point is that, in general, $c$ is not itself  the identity. However,
it is easy to see that it is described, up to homotopy, by the following composition
\begin{eqnarray}\label{equ:homeo}
W(Z)\times\Omega M\stackrel{\Delta\times id}{\longrightarrow}W(Z)\times
W(Z)\times\Omega M \stackrel{id\times p_{2}\times id}{\longrightarrow}\\
\to W(Z)\times Z(P,Q)\times\Omega M \stackrel{id\times
 j(P,Q)\times id}{\longrightarrow} W(Z)\times\Omega M\times\Omega
M\stackrel{id\times\nu}{\longrightarrow}  \nonumber \\  \to W(Z)\times\Omega M
\hspace{1.5in} \nonumber
\end{eqnarray}
where $p_{2}:W(Z)\to Z(P,Q)$ is the projection $D^{q}\times Z(P,Q)\to Z(P,Q)$
and, as before, $\nu$  is the multiplication. The restriction of $c$ to $\partial W(Z)$
is clearly also a homeomorphism and its description is obtained from (\ref{equ:homeo})
by restriction (recall that $\partial W(Z)=S^{q-1}\times Z(P,Q)\subset D^{q}\times
Z(P,Q)=W(Z)$). From (\ref{equ:homeo}) and (\ref{equ:tot_th}) we obtain the
following description of the homotopy class of $\delta'$
\begin{eqnarray}\label{equ:delt'}
S^{p-1}\times\Omega M\longrightarrow W(Z)\times\Omega M/\partial W(Z)\times \Omega
M\stackrel{\overline{c}}{\longrightarrow} \\ \to W(Z)\times\Omega M/\partial W(Z)\times
\Omega M \longrightarrow U(Q)\times\Omega M/A'(Q)\times\Omega M  \to \nonumber\\
\to U(Q)\times \Omega M/A(Q)\times\Omega M\stackrel{t}{\longrightarrow} S^{q}\wedge
\Omega M^{+}\nonumber
\end{eqnarray}
where $\overline{c}$ is the homeomorphism induced by $c$, $t$ is the obvious
composition $U(Q)\times\Omega M/A(Q)\times \Omega M\to (D^{q}\times D^{n-q}\times
\Omega M)/(S^{q-1}\times D^{n-q}\times\Omega M)\to (D^{q}\times \Omega
M)/(S^{q-1}\times\Omega M)\to S^{q}\wedge \Omega M^{+}$ and the map preceding $t$
is induced by a homeomorphism deforming $A'(Q)$ to $A(Q)$. We did also identify
$S^{u}(P)$ to $S^{p-1}$.

Let us denote by $\delta''$ the restriction of $\delta'$ in (\ref{equ:delt'}) to
$S^{p-1}\times{\ast}$. In view of the description of $c$ in (\ref{equ:homeo}), to show that
$\delta$ is the  composition in (\ref{equ:TP_comp}) it is enough to prove that $\delta''$ is
precisely the Thom map $h(P,Q)$ associated to the framing $\iota^{+}$ and the map $j(P,Q)$.
Recall that this Thom map is defined as the composition
\begin{equation}\label{equ:Th_m}
S^{p-1}\longrightarrow W(Z)/\partial
W(Z)\stackrel{\Delta}{\longrightarrow}\end{equation}
\vspace{0.05in}
$$\longrightarrow (W(Z)\times W(Z))/(\partial W(Z)\times W(Z))\stackrel{id\times
j(P,Q)}{\longrightarrow}$$
\vspace{0.05in}
$$\longrightarrow (W(Z)\times\Omega M)/(\partial W(Z)\times \Omega
M)\longrightarrow $$
\vspace{0.05in}
$$\stackrel{v}{\longrightarrow} (D^{q}\times Z(P,Q))\times\Omega M/(S^{q-1}\times
Z(P,Q))\times \Omega M)\stackrel{p_{1}\times id}{\longrightarrow} $$
\vspace{0.05in}
$$\longrightarrow D^{q}\times \Omega M/S^{q-1}\times \Omega M= S^{q}\wedge \Omega
M^{+}$$

\

where the map $v$ is the identification of $W(Z)$ with $D^{q}\times Z(P,Q)$
provided by the framing $\iota^{+}$. Notice that $\partial W(Z)=S^{q-1}\times Z(P,Q)$
is transported by the flow $\gamma$ into $\partial A(Q)\subset A(Q)$ in such a
way that $S^{q-1}\times \{x\}\subset \partial W(Z)$ is transported into $S^{q-1}\times
\{x'\}\subset S^{q-1}\times D^{n-q}=A(Q)$ by a homeomorphism.
This implies that the compositions of the last two maps in (\ref{equ:Th_m}) and the last three
in (\ref{equ:delt'}) are homotopic. In turn, this shows that $\delta''\he h(P,Q)$ and
concludes the proof of the theorem.

\

In the following we will consider Hopf complexes of certain particular pairs of spaces.
The first basic example is given by index pairs $(N_{1},N_{0})$ of some isolated invariant set
in a gradient flow on $M$. A cell decomposition of $N_{1}$
relative to $N_{0}$ produces as in  Corollary \ref{cor:hopf} a complex with coefficients
in $\pi_{\ast}^{S}(\Omega N_{1}^{+})$. We fix the following convention:
whenever talking about a Hopf complex of such an index pair we understand the complex
obtained after changing coefficients via the inclusion map $N_{1}\to M$. Of course,
this produces a complex with coefficients in $\pi_{\ast}^{S}(\Omega M^{+})$.
The second basic example of such pairs is given by the pair $(M',M)$ providing the
global Conley index of some isolated invariant set
 (see \S\ref{subsubsec:conlind}). Again, whenever talking about a Hopf complex of such a
global Conley index we understand the complex with coefficients in
$\pi_{\ast}^{S}(\Omega M^{+})$ obtained from a relative
cell decomposition of $(M',M)$  with the coefficients changed by the canonical projection map
$p:M'\to M$ of the global index.

\begin{corollary}\label{cor:ext_mo} The object $(C^{I}_{N}(f),d^{I}_{f})$ associated to a
local Morse function $(f,N)$ and introduced in Defintion \ref{extendedmorse} is indeed a chain
complex and is identified with a Hopf complex of an index pair for $I_{N}(f)$.
\end{corollary}

\begin{proof} Notice that in the setting of Definition \ref{extendedmorse}
if $x\in Crit^{N}_{a_{i}}(f)$ and $y\in Crit^{N}_{a_{i-1}}(f)$ for
$I=(a_{1},\ldots,a_{i-1},a_{i},\ldots )$, then $x$ and $y$ are consecutive critical points of
$f$ in $N$. Use the theorem to identify the bordism classes of the connecting manifolds of type
$Z(x,y)$ with the appropriate Hopf invariants and apply Corollary \ref{cor:hopf}.
\end{proof}

\begin{remark}\label{rem:ba} a. A different proof of $(d^{I}_{f})^{2}=0$ follows by adapting
to the non-compact case the method described in
\cite{Co2}. This is much more geometric in nature and is based on a precise description of the
boundaries of the spaces of smooth flow lines joining two critical points that are not
consecutive but might be joined by flow lines that are broken precisely once. One can also find
there an explicit discussion of the signs in the formula $(d^{I}_{f})^{2}=0$.

b. Clearly, Theorem \ref{theo:hopf_conn} used in conjunction with Proposition
\ref{prop:hopf} prove a bit more than the vanishing of the square of
the differential of $C^{I}_{N}(f)$. Indeed,
we obtain some relations among
the  Thom maps (of type $h(P,Q)$) associated to the connecting manifolds and not only
among their stabilizations. It is harder to encode these relations in a Morse type complex
the main problem being to find an appropriate coefficient ring (necessarily
non-commutative).  Moreover, even at the level of the bordism classes we get more
than the relations needed to show that $C^{I}_{N}(f,\alpha)$ is a chain complex.
Explicitely, if $P$ and $R$ are critical points that are connected by flow lines that
are broken at most once, then we have $\sum_{Q\in Crit_{\ast}(f)}+/-[Z(P,Q)][Z(Q,R]=0$ (the
precise signs are in \cite{Co2}; recall that $[Z(x,y)]$ is taken to be $0$ whenever
$x$, $y$ are not consecutive).

c. In view of the definition of the global Conley index
$\overline{c}_{N}(f)=\overline{c}_{\gamma(f)}(I(f))$ we also obtain that
$(C^{I}_{N}(f),d^{I}_{f})$ is identified with a Hopf complex of the global index because a
relative cell decomposition of an index pair for $I_{N}(f)$ produces also a relative cell
decomposition of the global Conley index of this invariant set.
\end{remark}

\subsection{Rigidity and duality for extended
Morse complexes.} \label{sec:rig}

\

In the first paragraph of this subsection we discuss certain continuation related properties of
extended Morse complexes. The second paragraph is concerned with the behaviour
of these complexes with respect to time reversion and stabilization.

\subsubsection{Continuation.}

\

Continuation induces a natural equivalence relation inside the space $\x(M)$
of pairs $(f,N)$ such that $f:M\to {\bf R}$ is a $C^{2}$ -function
and $N$ is an isolating neighbourhood for $\gamma(f)$  which, we recall, is the
flow induced by the negative gradient of $f$. As before, we denote by $I_{N}(f)$ the
maximal invariant set of this flow inside $N$. A pair $(f,N)\in\x(M)$
 will be called an {\em isolated local function on $M$}.

\begin{definition}\label{def:cont}
The isolated local functions  $ (f_{0},N_{0})$, $(f_{1},N_{1})$ are
 continuation equivalent if $I_{N_{0}}(f_{0})$ and $I_{N_{1}}(f_{1})$
viewed as isolated invariant sets in the flows $\gamma(f_{0})$ and respectively
$\gamma(f_{1})$ are related by continuation (see \S\ref{subsubsec:conlind}).
\end{definition}

\begin{remark}
a. Clearly, both the Conley index and the global Conley index are invariants
with respect to continuation equivalence.

b. Given $(f,N)\in \x(M)$ there is a $C^{2}$-neighbourhood $\mathcal{U}$ of
$f$ (in the weak topology) such that if $g\in \mathcal{U}$ then $(g,N)$
is also an isolated local function and is continuation equivalent to $(f,N)$.
\end{remark}

Here is a useful consequence of the results in
\S\ref{subsubsec:cor_hopf}

\begin{corollary}\label{theo:rig}
Assume that two local Morse functions $(f_{0},N_{0})$ and $(f_{1},N_{1})$
are continuation equivalent.
There exists a morphism of graded abelian groups
$m:{\bf Z}[Crit^{N_{0}}_{\ast}(f_{0})]\to {\bf Z}[Crit^{N_{1}}_{\ast}(f_{1})]$ wich
extends to a  chain-equivalence of free $\Omega^{fr}_{\ast}(\Omega M)$-modules
$$m^{I}:C^{I}_{N_{0}}(f_{0})\longrightarrow C^{I}_{N_{1}}(f_{1})$$
for each $I\subset {\bf N}$ for which both extended Morse complexes are defined.
\end{corollary}

\begin{proof}
Assume $C^{I}_{N_{0}}(f_{0})$ and $C^{I}_{N_{1}}(f_{1})$ are both defined.
Because $(f_{0},N_{0})$ and $(f_{1},N_{1})$ are continuation equivalent
then the global Conley indexes (see \S\ref{subsubsec:conlind} )
$\overline{c}_{N_{0}}(f_{0})$ and $\overline{c}_{N_{1}}(f_{1})$ coincide.
The two extended Morse complexes above being identified by  Corollary \ref{cor:ext_mo}
to Hopf complexes of this same global index, the statement follows by
applying Lemma \ref{lem:h_comp}.
\end{proof}

\begin{remark} There is also a purely geometric proof of a variant of this corollary  which
 is based on the techniques that appear in \cite{CoRa}.
The geometric method leads in fact to more powerful results because it may be
applied in cases when the global Conley index is not defined. Moreover,
it also shows that  for any Morse function $(f,N)$
an arbitrary perturbation of the metric $\alpha$
inside $N$ does not modify the isomorphism type of the extended Morse complexes of $f$.
\end{remark}

We will make use of the following notion of minimality.

\begin{definition}\label{def:min} A local Morse function $(f,N)$ is
minimal if $C^{\bf N}_{N}(f)$  has trivial differential. In other
words, the classical Morse complex of $f|_{N}$  has
trivial differential.
\end{definition}

If  a local Morse function $(f,N)$ is minimal, then its number of critical points $\in N$ of
index $i$ is equal to $rk(H_{i}(c_{N}(f);{\bf Z})$. Obviously,  the integral homology of
$c_{N}(f)$ is forced to be torsion free. The next statement immediately follows from
 Corollary \ref{theo:rig}.

\begin{corollary}\label{cor:PS_min} Assume that two local Morse functions $(f,N)$, $(g,N')$
are minimal and continuation equivalent.
Then for any index set $I$ we have $C_{N}^{I}(f)\approx C_{N'}^{I}(g)$.
\end{corollary}

At this point we introduce a particular notation for certain special extended Morse complexes.
Let $(f,N)$  be a local Morse function and let $I_{f}=\{k\in {\bf N} :
Crit^{N}_{k}(f)\not=\emptyset\}$. We denote by $(\mathcal{C}_{N}(f),d_{f})$ the
extended Morse complex  $(C^{I_{f}}_{N},d^{I}_{f})$. This complex carries the most
possible information in its differential among all extended Morse complexes of $(f,N)$.
Recall the minimal Hopf complexes $\mathcal{C}_{min}(X',X)$ that appear
in Lemma \ref{lem:h_comp} and are defined whenever  $(X',X)$ admits a minimal
cell-decomposition in the sense of that lemma.

From the results above we immediately obtain.

\begin{corollary}\label{cor:min}
Assume that the local Morse function $(f,N)$ is  minimal.
Then $\mathcal{C}_{min}(\overline{c}_{N}(f))$
is defined and is isomorphic to $\mathcal{C}_{N}(f)$.
\end{corollary}

\begin{remark} a. In practice, to apply the statement above one needs to
construct a minimal cell decomposition for the global Conley index $\overline{c}_{N}(f)$
purely homotopically (this is possible, for example, if this Conley index has torsion
free integral homology and is simply connected) and then compute
$\mathcal{C}_{min}(\overline{c}_{N}(f))$ out of it.

b. One can, of course, use instead of the global Conley index an index pair
$(N_{1},N_{0})$ for $I_{N}(f)$.
\end{remark}

The existence of a minimal Morse function in a given continuation class
 can be insured only under circumstances somewhat restrictive.
The result below follows, as shown in \cite{Re},
from the classical techniques used by Smale in the proof of the h-cobordism theorem.

\begin{proposition}\cite{Re} Assume the isolated local function $(f,N)$ has the
following two properties:
\begin{itemize}
\item[(i)] $c(f)$ has torsion free integral homology
\item[(ii)] there exists a bi-regular index block
$(N_{1},N_{0}, N_{2})$ for $I_{N}(f)$  with  $N_{i}$  simply connected for $i\in\{0,1,2\}$.
\end{itemize}
Then there is a minimal  local Morse function continuation equivalent to $f$.
\end{proposition}

\subsubsection{Stabilization and duality of extended Morse complexes.}\label{subsub:du_m}

\

For $(f,N)\in\x(M)$, $(h,N')\in\x(M')$ we may define $f\oplus h:M\times M'\to {\bf R}$ by
$(f\oplus h)(x,y)=f(x)+h(y)$ and clearly, $(f\oplus h, N\times N')\in\x(M\times M')$ (where the
metric on the product is the product metric). Moreover, by the basic properties of the Conley
index we have $c_{N\times N'}(f\oplus h)=c_{N}(f)\wedge c_{N'}(h)$. Of particular iterest is
the case when $h$ is a non-degenerate quadratic form $q:{\bf R}^{k}\to {\bf R}$ the metric on
${\bf R}^{k}$ is the euclidean one and $N'$ is a disk $D^{k}$ containing $0$. A function of
the form $f\oplus q$ is called a {\em stabilization} of $f$ \cite{ElGr}\cite{Co1}\cite{Co2}.
 Of course, in this case  $c_{N'}(q)=S^{ind(q)}$.

The behaviour of the extended Morse complexes with respect to stabilization is
very easy to describe.

\begin{corollary}\label{cor:susp} Let $(f,N)$ be a local Morse function.  Assume
$q:{\bf R}^{k}\to {\bf R}$  is a non-degenerate quadratic form.
For any index set $I\subset {\bf N}$ we have an isomorphism
$s^{ind(q)}(\mathcal{C}^{I}_{N}(f),d^{I}_{f})\approx \mathcal{C}_{N\times
D^{k}}^{I'}(f\oplus q,), d^{I'}_{f\oplus q})$ where $I'=\{k+ind(q) : k\in I\}$ and $s^{(-)}$ is
algebraic suspension.
\end{corollary}

\begin{proof}  The statement is immediate by noting that
$f+q$ is Morse-Smale with respect to the sum metric $\alpha+euclidean$ on
$N\times  D^{k}$. Moreover, two consecutive
critical points $P$ and $Q$ of $f$ remain consecutive for $f\oplus q$ and have the same
connecting manifold. The framing of this connecting manifold also does not change stably and
the map to the loop space $\Omega(M\times {\bf R}^{k})\simeq \Omega M$ also remains the
same up to homotopy.
\end{proof}

It was shown in the first paper in this series \cite{Co1} that if the stable normal
bundle of the manifold $M$ (assumed there to be compact)  is trivial, then for any flow
$\gamma$ on $M$ and any  isolated invariant set $S$ of $\gamma$ the Conley indexes of $S$
with respect to the direct flow $c_{\gamma}(S)$ and with respect to the inverse flow
$c_{-\gamma}(S)$ are related by Spanier-Whitehead duality (here
$(-\gamma)_{t}(x)=\gamma_{-t}(x)$). This duality relation was shown to be also valid, in this
situation, for the homotopical data resulting from an attractor-repellor pair. In particular the
connection maps of the attractor and the repellor with respect to the direct and inverse flows are
Spanier-Whitehead duals. A partial, homotopical variant of this result in the Morse case was
already present in \cite{Fr}. A homological version  in the case of general isolated invariant
sets appeared in \cite{Mc}. If the stable normal bundle of
$M$ is no longer trivial, then it is natural to ``estimate" the possible lack of duality. This
was investigated in \cite{Co1} by means of a certain Thom construction.

The last task  of this sub-section
 is to use our understanding of connecting manifolds to explicit the deviation
from Spanier-Whitehead duality (along
the lines in \cite{Co2}) for the case of the simplest type of attractor-repellor pair, consecutive
non-degenerate critical points in a Morse-Smale gradient flow, and to then formulate the result
globaly in terms of the extended Morse complexes.

On the purely homotopical side some duality results already appear in
Corollary \ref{cor:PD} . Here we treat duality from the perspective of connecting manifolds
and extended Morse complexes.

We return now to our manifold $M$ and we asssume further in this paragraph that there is a
compactification $\overline{M}^{n}$ of
$M$ which is included in a compact, closed manifold $M'$ of dimension $n$. Let $\nu$ be the
stable normal bundle of $M'$ restricted to $M$. We will also denote by $\nu$ the resulting
map $\nu:M\to BSO$.  We will refer to this bundle $\nu$ as a stable normal bundle to $M$.

There is an
involution $\ast : \pi_{\ast}^{S}(\Omega M^{+})\to \pi_{\ast}^{S}(\Omega M^{+})$ defined
as follows: if $a\in \Omega^{fr}_{k}(\Omega M)$ is given by a map $a: S^{t}\to
S^{t-k}\wedge \Omega M^{+}$, $t>>0$, then $a^{\ast}$ is represented by the composition
$S^{t}\stackrel{a}{\to}S^{t-k}\wedge \Omega M^{+}\stackrel{id\wedge
\Delta}{\longrightarrow}
S^{t-k}\wedge \Omega M^{+}\wedge \Omega M^{+}\stackrel{id\wedge
\Omega\nu \wedge id}{\longrightarrow} S^{t-k}\wedge {\bf SO}^{+}\wedge\Omega M^{+}
\stackrel{\mu\wedge id}{\longrightarrow} S^{t-k}\wedge \Omega M^{+}\stackrel{id\wedge
-1}{\longrightarrow} S^{t-k}\wedge\Omega M^{+}$ where $\mu$ is induced by the action of
${\bf SO}(t-k+1)$ on $S^{t-k}$,  $-1$ is the map reversing the orientation of the loops in
$\Omega M$ and $\Delta$ is the diagonal.

Assume that $(\mathcal{C},d)$ is an $\pi_{\ast}^{S}(\Omega M^{+})$-chain complex. We
define its $n$-dual $(\mathcal{C}^{n},d^{\ast})$ by an isomorphism
$\ast:\mathcal{C}^{\ast}_{n-k}\approx \mathcal{C}_{k}$ and
$d^{\ast}(x^{\ast})=\sum d^{\ast}(x^{\ast},y^{\ast})y^{\ast}$ such that
$d^{\ast}(x^{\ast},y^{\ast})=(d(y,x))^{\ast}$ with $d(y)=\sum d(y,x)x$.

\begin{corollary} \label{cor:du_ext} Assume $(f,N)$ is a local Morse function on $M$.
\begin{itemize}
\item[(i)] Suppose that $P,Q$ are consecutive critical points of $f|_{N}$. If
$[Z(Q,P)]\in\Omega^{fr}_{\ast}(\Omega M)$ is the bordism class of the connecting manifold
of $Q$ and $P$ viewed as consecutive critical points for the function $-f$, then
$[Z(Q,P)]=[Z(P,Q)]^{\ast}$.
\item[(ii)] we have an isomorphism $(\mathcal{C}_{N}(-f),d_{-f})\approx
(\mathcal{C}_{N}(f)^{\ast},(d_{f})^{\ast})$.
\item[(iii)] if $f$ is a minimal Morse function and $M$ is closed, $N=M$, then we have
$\mathcal{C}_{N}(f)\approx (\mathcal{C}_{N}(f))^{\ast}$ the isomorphism being induced
by the Poincar\'e duality isomorphism.
\end{itemize}
\end{corollary}
\begin{proof} The first point is an easy consequence of the fact that, as a toplogical space, the
connecting manifold of $P$ and $Q$ viewed as consecutive critical points for $f$ is the same as
that of $Q$ and $P$ as consecutive critical points for $-f$. Moreover, the inclusion into
$\Omega M$ changes by the reversion of the loops. It remains to compare the
framings of $Z(P,Q)$ in $S^{u}_{f}(P)$ and $S^{s}_{f}(Q)$ (the subscript $f$ indicates
that the flow used is $\gamma(f)$ ). Let $E(M')$ be a tubular neighbourhood of $M'$
inside a high dimensional sphere. Of course, $E(M')$ is diffeomorphic to the total space
of the disk bundle associated to $\nu$. Let $q:E(M')\to {\bf R}$ be a function measuring
the square of the distance from $M'$ inside $E(M')$. Inside each fibre of the disk
 bundle $E(M')$ the function $q$ is a positive definite quadratic form.
 Let $E(M)$ be the restriction of the bundle
to $M$ and  consider now the function $\overline{f}: E(M)\to {\bf R}$, $\overline{f}(z)=
f(p(z))+q(z)$ where $p:E(M')\to M'$ is the projection of the bundle. It is easy to see that
$\overline{f}$ is a local Morse function after choosing an appropriate metric on $E(M)$.
Assume that the rank of $E(M)$ is $k$. Notice that
$S^{u}_{\overline{f}}(P)=S^{u}_{f}(P)$ and
$S^{s}_{\overline{f}}(Q)=\Sigma^{k}S^{s}_{f}(Q)$. Because a sphere has a trivial normal
bundle the framings of $Z(P,Q)$ in $S^{u}_{\overline{f}}(P)$ and
$S^{s}_{\overline{f}}(Q)$ are stably the same up to sign (by the more general results in
\cite{Co1} or those in \cite{Fr}). On the other hand the framings of $Z(P,Q)$ in
$S^{u}_{f}(P)$ and $S^{u}_{\overline{f}}(P)$ are the same. However, the framings
of $Z(P,Q)$ in $S^{s}_{f}(Q)$ and $S^{s}_{\overline{f}}(Q)$ differ by the twisting
coming from the normal bundle of $M'$ in $E(M')$ which is $\nu$.  This immediately
implies (i). The second point is an obvious
consequence of the first. The third is implied by the fact that both $f$ and $-f$ are
minimal and as $M$ is compact they are in the same continuation class.
\end{proof}

\begin{remark}  Notice that, in the smooth context, this corollary is an extension of Corollary
\ref{cor:PD}.  It can be used to detect homotopy equivalences that do not admit
diffeomorphisms as representatives. The third point can be interpreted as saying that the
Poincar\'e duality isomorphism respects the structure contained in the minimal Hopf
complex of $M$.
\end{remark}

\section{Bounded, closed and periodic orbits for hamiltonian flows.}\label{sec:orb}

\

The setting is as follows. We assume from now on that our
fixed manifold $M$ admits a symplectic form $\omega$ that we fix.
 Assume that $f:M\to {\bf R}$
is a smooth function. The hamiltonian vector field induced by $f$, $H_{f}$, is defined
by requiring the equation $\omega(H_{f},Y)=-df(Y)$ to hold for all vector fields $Y$ on
$M$. The flow induced by the vector field $H_{f}$ will be denoted by $h^{f}$ and is
called the hamiltonian flow of $f$.

In searching for bounded and  periodic orbits the following immediate remark is useful.
For any riemannian metric $\alpha$ on $M$
we have $\alpha (-\nabla^{\alpha}(f),H_{f})=-df(H_{f})=\omega(H_{f},H_{f})=0$ thus
 $H_{f}$ and $\nabla^{\alpha}(f)$ are $\alpha$-orthogonal and each hypersurface
$f^{-1}(y)$ is $h^{f}$ -invariant. Moreover, it is well known
that the problem of existence of periodic orbits of $h^{f}$ on some regular hypersurface
$A=f^{-1}(a)$ does only depend on $A$ and not on $f$. This is why one can expect that some
topological constraints on $f$ can impose enough tension on the hypersurfaces $f^{-1}(a)$ as
to force some recurrence phenomena.

The case that has been extensively studied in the literature corresponds to a
function $f$ that has some compact regular hypersurfaces. The simplest way to insure this
condition is to take $f$ to be a function that attains an  extremum on a
compact set.   Suppose $A=f^{-1}(a)$ is a regular, compact hypersurface. Then, of course,
every orbit of $h^{f}$ originating in $A$ is bounded, Poincar\'e's recurrence theorem
shows that, except for a set of zero measure, all points are recurrent. As for the existence of
periodic orbits there are two complementary approaches.
The first, based on the $C^{1}$ closing
lemma of Pugh and Robinson, shows that for a generic choice of $f$ the recurrent points can be
transformed into periodic ones and thus, in this compact case, ``generically", periodic orbits are
abundant.  The second approach has been pursued succesfully by Rabinowitz
\cite{Ra}, Weinstein \cite{We} as well as Ekeland, Floer, Hofer, Moser, Zehnder and
many other authors (see \cite{HoZe} and \cite{Ho} for surveys on the subject).
It consists in using methods adapted to the specific form of a fixed, particular $f$ and
namely variational methods on the space of  the free loops in $M$ to
deduce the existence of  periodic orbits for this fixed $f$.

\

In the following we focus on the non-compact context. We assume therefore
that $M$ is {\em non-compact} and, as indicated in the introduction, at the heart
of our results are certain criteria that will insure the existence of bounded orbits.
In the first sub-section it is the lack of (Spanier-Whitehead) duality which leads to such
existence results. In the second sub-section we use the non-vanishing
of the differential in an extended Morse complex, or (given the results
in Section \ref{sec:conn}) equivalentely, the presence of a non-vanishing relative Hopf invariant.
It should be pointed out that the key advantage of
the second method is that the bounded orbits that are produced in this way are better
localized, in particular, they belong to regular hypersurfaces, fact that can not be
guaranteed by the lack-of-duality approach.
We then use the $C^{1}$-closing lemma to obtain {\em generically} periodic or,
(in the case of the  lack of Spanier-Whitehead  duality method) closed (possibly
homoclinic) orbits.  Compared to the understanding of the compact case our results only
correspond to the first, ``generic'', point of view.
By following the arguments below it is easy to see that, in fact, both
methods used to produce bounded orbits  actually apply not only to
hamiltonian flows but to any flows that are orthogonal to gradient ones.
The third sub-section contains some explicit examples as well as the
interpretation of our various results in terms of the invariant $d(-)$ mentioned in the
introduction.

It is to be  expected that by adding some
appropriate analytical assumptions to our homotopical  conditions one will obtain non-generic
existence results. For example, I do not know of any examples of contact hypersurfaces
that carry bounded characteristics and do not carry closed ones.

\

For now, let us  notice that the bounded orbit existence problem is non-trivial.

\begin{example}\label{exam:non_bded}
Let $M={\bf R}^{2k}$ be endowed with the
symplectic form $\omega_{0}=\sum_{i=1}^{k}dx_{i}\wedge dy_{i}$ and let
$f(x_{1},y_{1},\ldots ,x_{k},y_{k})= \sum_{i=1}^{k}f_{i}(x_{i},y_{i})$
with $f_{i}:{\bf R}^{2}\to {\bf R}$ a smooth, analytic function with a single
singularity at $0\in{\bf R}^{2}$ which is different from an extremum, $1\leq i\leq k$.
Then $h^{f}$ has no bounded, non-trivial orbits. Moreover, the family of functions $f':{\bf
R}^{2k}\to {\bf R}$ such that $h^{f'}$ has some non-trivial bounded orbit is not dense in any
$C^{2}$-(weak) neighbourhood of $f$.  It is useful to remark at this point that, if $0$ is
a totally degenerate critical point of $f$, then arbitrarily close to $f$ ( in the strong
$C^{2}$ topology)
one can find  functions $g$ that have some compact level surfaces (one can construct such a
function by just forcing the appearance of a local maximum close to the origin in ${\bf R}^{2k}$)
and thus the hamiltonians of type $h^{g}$ will have bounded non-trivial orbits.
It is easy to see that in the example above, under the additional assumption that
$0$ is a non-degenerate critical point, the linearized hamlitonian flow at $0$ is defined
and it does not have any (non-trivial) periodic orbits. However, there are more
sophisticated examples (the most famous produced by Moser \cite{Mo})
where no bounded orbits of the hamiltonian flow
exist even if the linearized flow at the singular point has  plenty of periodic ones.
\end{example}

\subsection{Lack of self - duality and closed orbits.}

\

Recall that for $(f,N)\in\x(M)$ we denote by $c_{N}(f)$ the Conley index of the maximal
invariant set  $I_{N}(f)$ of $f$; $\overline{c}_{N}(f)$ is the global Conley index (both with
respect to the negative gradient of $f$).

We start with a result that continues an idea that has first been exploited in \cite{Co3}.
Recall from sub-section \ref{subsub:du_m} that a stable normal bundle of $M$
is the stable normal bundle of a closed manifold $M'$ of same dimension as $M$
and which contains a compactification of $M$. Using this normal bundle recall also
that we can construct the dual of an extended Morse complex .

\begin{proposition} \label{prop:bded_orbb} Assume that a normal bundle $\nu$ to $M$ is
defined. Suppose that $(f,N)\in\x(M)$ is fixed and suppose that it has the property that
$I_{N}(f)\subset f^{-1}(v)$ for some isolated critical value $v\in{\bf R}$ of $f$.
Suppose that one of the following conditions is satisfied:
\begin{itemize}
\item[(i)]  the restriction of $\nu$ to an isolating neighbourhood of $I(f)$ is trivial
and $c_{N}(f)$ is not Spanier-Whitehead self $n$-dual.
\item[(ii)]  $f$ is continuation equivalent to a minimal Morse function,\\
and $\mathcal{C}_{min}(\overline{c}_{N}(f))$  is not  isomorphic to
$(\mathcal{C}_{min}(\overline{c}_{N}(f))^{\ast}$.
\end{itemize}
Then the hamiltonian flow $h^{f}$ has the property that for any
isolating neighbourhood $K$ of $I_{N}(f)$ there is a bounded orbit of $h^{f}$ that is
included in $K$ and intersects the boundary of $K$.
\end{proposition}
\begin{proof}
We first claim that each of the conditions at (i) and (ii) implies that the invariant sets
$I_{N}(f)$ and, respectively, $I_{N}(-f)$  are not related by continuation.
Indeed, if they would be related by continuation, then
$c_{N}(f)=c_{N}(-f)$. But as $\nu$ is trivial on an isolating neighbourhood of the maximal
invariant set of $f$, by results in \cite{Co1}, we have that $c_{N}(-f)$  is an $n$-
Spanier-Whitehead dual of
$c_{N}(f)$ hence, our claim follows from point (i) (we recall that the spaces $A$,
$B$ with the homotopy type of finite $CW$-complexes are $n$-Spanier-Whitehead duals if
there exist $A',B'$ complementary in an $S^{2m+n+1}$ sphere and with $A'\simeq
\Sigma^{m}A$, $B'\simeq \Sigma^{m}B$).
Similarly, condition (ii) also implies
the non-existence of this continuation by applying Corollary \ref{cor:du_ext} to the minimal
Morse function continuation equivalent to $f$ and equating its minimal extended Morse
complex with $\mathcal{C}_{min}(\overline{c}_{N}(f))$ (by Corollary
\ref{cor:min}). This Hopf complex only depends on the global Conley index of $f$. We
then notice that if $I_{N}(f)$ and $I_{N}(-f)$ are related by continuation, then their global
Conley indexes are the same.  Applying the same argument to the negative of the minimal
Morse function in the continuation class of $f$ we obtain our claim.

We now arrive to the key idea of the proof.  We  notice that if the conclusion of
the proposition would be false, then there would exist  a continuation between
$I_{N}(f)$ and $I_{N}(-f)$.  Indeed, we let $u,v,w:[-1,1]\to [0,1]$ be a smooth partition of the
unity such that $u^{-1}(0)=[0,1]$, $v^{-1}(0)=[-1,0]$, $w^{-1}(0)=[1,-1/2]\cup [1/2,1]$
$u([-1,-1/2])=1$, $v([1/2,1])=1$, $w([-1/4,1/4])=1$, $u'(t)<0$ for $t\in (-1/2,0)$,
$v'(t)>0$ for $t\in (0,1/2)$, $w'(t)sign(t)<0$ for $t\in (-1/2,-1/4)\cup (1/4,1/2)$.
We now define a smooth one parameter family of vector fields on $M$, $X_{\tau}=
-u(\tau)\nabla (f)+w(\tau)H_{f}+v(\tau)\nabla (f)$.  Suppose that $K$ is
an isolating neighbourhood of $I_{N}(f)$ which does not satisfy the conclusion of the
proposition. Then $K\times [-1,1]$ is an isolating neighbourhood for the flow
$\Gamma$ on $M\times [-1,1]$ given by letting $\Gamma_{\tau}=\Gamma|_{M\times
\{\tau\}}$ be induced by $X_{\tau}$. This happens because for all values of $\tau\not=0$
one of the functions $f$ or $-f$ is a Lyapounov function for the flow $X_{\tau}$ and as
$I_{N}(f)$ is contained in a single level hypersurface of $f$ it results that the maximal compact
invariant set for $\Gamma_{\tau}$ inside $K\times\{\tau\}$ is $I_{N}(f)$. For $\tau=0$ our
assumption that all orbits of $h^{f}$ that are contained in $K$ do not intersect $\partial K$
implies that $K\times\{0\}$ is an isolating neighbourhood for $\Gamma_{0}$.
This means that $I_{N}(f)$ can be continued to $I_{N}(-f)$ and leads to a contradiction.
\end{proof}

By a closed orbit in a flow we mean a periodic orbit or the trajectory of some
point $x$ whose $\omega$-limits $\omega(x)$ and $\omega^{\ast}(x)$ both
reduce to the same single point (this last case is refered to as a homoclinic closed orbit).

\begin{corollary} \label{cor:hom} Let $(f,N)\in\x(M)$. There is a family of functions
$\mathcal{F}$ dense in a $C^{2}$- neighbourhood of $f$  in the  $C^{2}$-strong topology
such that for each
$f'\in\mathcal{F}$ the hamiltonian flow $h^{f'}$ has at least
$$\sum_{j\not=n/2}rk (H_{j}(c(a); {\bf Z})$$ possibly homoclinic, non-trivial
closed orbits.
\end{corollary}
\begin{proof}
We fix a sufficientely small $C^{2}$ -neighbourhood $\mathcal{U}$ of $f$ such that
any function $g\in \mathcal{U}$ has the property that $N$ is an isolating neighbourhood
for $\gamma(g)$. There is dense family $\mathcal{F}'\subset\mathcal{U}$  such that
$\mathcal{F}'$ is dense in $\mathcal{U}$, each critical level of a
function $f'\in\mathcal{F}'$ only contains a single critical point inside $N$, $(f',N)$ is a local
Morse function and is continuation equivalent to $f$ (this is true because starting with any
isolated local function $(g,N)$ we may perform a compactly supported, arbitrarily small
modification to obtain a function
$f'$ as above). Let $f'\in\mathcal{F}'$. An immediate application of
the Morse inequalities shows that such a function has at least $rk(H_{k}(c_{N}(f);{\bf Z}))$
critical points of index $k$. We may apply Proposition \ref{prop:bded_orbb}  for each critical
point $P$ of  $f'$ such that $ind(P)\not=n/2$. We obtain that for each sufficiently small
compact set $K_{P}$ around $P$ the hamiltonian flow $h^{f'}$ has at least one non-trivial
bounded orbit in $K_{P}$ that intersects non-trivially the boundary of $K_{P}$.
This means that either this orbit has both its $\omega$ and $\omega^{\ast}$
limits equal to $P$ - in which case it is a homoclinic orbit - or one of these limits belongs
to $\Omega_{c}(h^{f'})\cap K_{P}-\{P\}$ where for a flow $\gamma$,
$\Omega_{c}(\gamma)$ is the set of
the non-wandering points of $\gamma$ that have a non-vanishing $\omega$ or
$\omega^{\ast}$ limit.  In this case, let $y\in \Omega_{c}(h^{f'})\cap
K_{P}-\{P\}$.  By  the $C^{1}$ closing lemma of Pugh and
Robinson \cite{PuRo} (in particular the argument in \S 11) there is a
$C^{1}$-deformation of $H_{f'}$ that has a compact support (only depending on the orbit of
$y$), that can be assumed arbitrarily small and whose result is a
hamiltonian vector field $H_{f''}$ whose associated flow has a periodic orbit that is close to
$y$ and therefore is non-constant. Because the deformation leading to $H_{f''}$ is compactly
supported it results that we may obtain $f''$ as $C^{2}$-close to $f'$ as desired and by apply
this construction for each such critical point $P$. This proves the existence of the family
$\mathcal{F}$.
\end{proof}

\subsection{Non-vanishing of Hopf invariants and periodic orbits.}

\

In this sub-section we search for bounded and periodic
orbits that are situated on regular level surfaces.

Let $(f,N)$ be a local Morse function such that $I_{N}(f)$ satisfies the property
$(\ast)$ from the introduction. Before stating the main theorem of this section we indicate the
role of property $(\ast)$ in this context. By definition, this property provides an index pair
$(N_{1},N_{0})$ of
$I_{N}(f)$ with $N_{0}\subset f^{-1}(a)$ with $a$ a  regular value of $f|_{N_{1}}$.
By using an appropriate Lyapounov function we see
that we may immediately assume that there is a strong index block $(N_{1},N_{0},N_{2})$
with $N_{0}\subset f^{-1}(a)$  that we fix.
For each $P\in Crit^{N}(f)$ we see that $W^{u}_{N}(P)\cap\partial
N_{1}=W^{u}_{N}(P)\cap\partial f^{-1}(a)$.  We will discuss criteria insuring that property
$(\ast)$ is satisfied in \S\ref{subsec:appl}.

Let $P,Q\in Crit^{N}(f)$ such that their indexes $p,q$  are
succesive in
$\{k: Crit^{N}_{k}(f)\not=\emptyset\}$. As in the setting of Theorem \ref{theo:hopf_conn}
we consider the relative Hopf invariant $H(j_{P},j_{Q})$. Recall that
we denote by $[H(j_{P},j_{Q})]\in \pi_{\ast}^{S}(\Omega M^{+})$ the stable
class of $H(j_{P},j_{Q})$.

\begin{theorem} \label{theo:closed} With the assumptions above, if
$[H(j_{P},j_{Q})]\not=0$, then  arbitrarily close to
$f$ in the strong $C^{2}$ topology there is a function $\overline{f}$ such that there are infinitely
many regular values $v$ of $\overline{f}$, $f(Q)<v<f(P)$ with the property that
$\overline{f}^{-1}(v)\cap N_{1}$ contains a periodic orbit of $h^{\overline f}$ .
\end{theorem}

\begin{remark} Because $f|_{N}$ is Morse, any function $f'$ sufficiently $C^{2}$-close to $f$
is also Morse after restriction to $N$ and has the same type and number of critical points as
$f$ there.
\end{remark}

The proof of the theorem  occupies the rest of this subsection

\

Assume the setting of  Theorem \ref{theo:closed}.
The main step is again to detect bounded orbits and consists in proving
the next statement.  In its turn, this has considerable
intrinseque interest as it applies to the
function $f$ itself (whithout requiring any perturbation). Moreover,
there exist  examples when $[H(j_{P},j_{Q})]=0$ but $H(j_{P},j_{Q})\not=0$.

\begin{proposition}\label{prop:bded_gen}
Assume that one of the following two conditions is satisfied
\begin{itemize}
\item[(i)] $H(j_{P},j_{Q})\not=0$ and there are no critical points $R\in Crit^{N}_{q}(f)$ with
$f(Q)<f(R)<f(P)$, $Z(P,R)\not=\emptyset$.
\item[(ii)] $[H(j_{P},j_{Q}]\not=0$ and  each $R\in Crit^{N}_{q}(f)$
with $f(Q)<f(R)<f(P)$ verifies $[H(j_{P},j_{R})]=0$.
\end{itemize}
then there exists $\delta>0$ such that for all $v\in (f(Q), f(Q)+\delta)$ the intersection
$V_{v}=f^{-1}(v)\cap Int(N_{1})$ is  regular and contains at least one bounded orbit of the
hamiltonian flow $h^{f}$ .
\end{proposition}

\begin{proof}

We start with a general simple, auxiliary result.

\begin{lemma}\label{lem:large_ind} Assume $\gamma:M'\times {\bf R}\longrightarrow M'$
is a flow that has an isolated invariant set $S\subset M$.
Assume $K\subset M'$ is an isolating neighbourhood of $S$. Let $K'$ be a compact
set such that $K\subset Int(K')$ and $I_{\gamma}(K')=S$. There is an index pair
$(L_{1},L_{0})$ for $S$ such that $K\subset Int(L_{1}-L_{0})$ and $L_{1}\subset K'$.
\end{lemma}
\begin{proof} This is a simple variation on the construction of index pairs in \cite{Sa}.
For a compact set $T\subset K'$ with $S\subset Int(T)$  let
$W_{T}=\{x\in K': \exists t\in {\bf R}^{+}\  , \  \gamma_{[0,t]}(x)\subset K',
\gamma_{t}(x)\in T\}$.
This set is compact. Let $K''\subset Int(K')$ be a second compact set with $K\subset Int(K'')$
and let $U\subset W_{K''}$ be an open neighbourhood of $W_{K}$. Let
$L_{0}=\{x\in K': \exists \ t\in {\bf R}^{+}, \gamma_{[-t,0]}(x)\subset K', \gamma_{-t}(x)\in
K'-U\}$. This set is positively invariant in $K'$, compact and $K\cap L_{0}=\emptyset$.
Let $L'_{1}=\{x\in K':\exists \ t\in {\bf R}^{+}, \gamma_{[-t,0]}(x)\subset K',
\gamma_{-t}(x)\in K'' \}$. This is again
compact and we take $L_{1}=L_{0}\cup L'_{1}$. It is easy to verify that $(L_{1},L_{0})$ is
an index pair.
\end{proof}

We return now to the proof of Proposition \ref{prop:bded_gen}.

To simplify notation we let $W^{u}(-)=W^{u}_{N_{1}}(-)$ and similarly for the stable
manifolds. In the same way we drop the index $N_{1}$ (or $N$) in the notation for the
various isolated invariant sets as long as everything takes place in the isolating neighbourhood
$N_{1}$. Moreover, because all our arguments will take place in $N_{1}$ we may
assume (after possibly modifying the vector field $H_{f}$ away from $N_{1}$) that
$h^{f}$ is a flow on $M$ (and not only a partially defined one).

The first step is to assume that
we are in a particular situation.

\begin{lemma}\label{lem:bded_ho}  If $H(j_{P},j_{Q})\not=0$ and all values $w\in
(f(Q),f(P))$ are regular for $f|_{Int(N_{1})}$, then for some $\delta'>0$  and  each $v\in
(f(Q),f(Q)+\delta')$ the regular hypersurface $V_{v}$ contains at least one
bounded orbit of the hamiltonian flow  $h^{f}$.
\end{lemma}

\begin{proof}
For a value $v\in (f(Q),f(P))$,  recall $V_{v}=f^{-1}(v)\cap Int(N_{1})$ be the corresponding
level hypersurface of $f$. The
condition $(\ast)$ implies that flow lines of $\gamma=\gamma(f)$ originating in $P$
intersect $V_{v}$. We take $\delta'$ sufficiently small such that for $v\in (f(Q),f(Q)+\delta')$
all the flow lines arriving in $Q$ also have to cross $V_{v}$. We now fix such a $v\in
(f(Q),f(Q)+\delta')$ and let $V=V_{v}$.
Let $S(P)=S^{u}(P)=W^{u}(P)\cap V$ be the unstable sphere
of $P$ and $S(Q)=S^{s}(Q)=W^{s}(Q)\cap V_{v}$  be the stable sphere of $Q$.
As  there are no critical values in the interval $(f(Q),f(P))$  we obtain
that $S(P)\approx S^{p-1}$ and $S(Q)\approx S^{n-q-1}$. Let $K'\subset V$ be a
compact neighbourhood of the union $S(P)\cup S(Q)$.

We intend to show by contradiction that $K'$ contains (at least) one orbit of $h^{f}$.
We now suppose that this is not the case and therefore the maximal compact invariant
set of $h^{f}$ inside $K'$, $I_{h^{f}}(K')$, is void.
Let $Z=Z(P,Q)=S(P)\cap S(Q)$. Recall that,  as in Theorem \ref{theo:hopf_conn}, we have
the framed embedding $i:Z\subset S(P)$ and the map $j=j(P,Q):Z(P,Q)\to \Omega M$.
Our purpose is to see that the assumption that $I_{h^{f}}(K')=\emptyset$ is sufficient to
construct a null-bordism of the the pair $(i,j)$ inside $S^{p-1}\times [0,1]$.
This immediately leads to a contradiction. Indeed,  if the couple $(i,j)$ is bordant to $0$
inside $S^{p-1}\times [0,1]$  then the Thom map associated to $(i,j)$,
$t=h(P,Q):S^{p-1}\to S^{q}\wedge (\Omega M^{+})$, is null-homotopic. But
Theorem \ref{theo:hopf_conn} claims that $h(P,Q)$ is homotopic (up to sign) to
$H(j_{P},j_{Q})$ which is non-trivial by hypothesis.

We notice that $K'$  is an isolating
neighbourhood and as $S(P)\cup S(Q)$ is compact we may find (by Lemma
\ref{lem:large_ind}) an index pair $(L_{1},L_{0})$  of  $I_{h^{f}}(K')$ such that $S(P)\cup
S(Q)\subset Int(L_{1}-L_{0})$, $L_{1}\subset K'$  (the fact that
$I_{h^{f}}(K')=\emptyset$  does not prevent the application of that lemma). Moreover, by the
construction of regular index pairs in \cite{Sa} we may even assume that $(L_{1},L_{0})$ is
a regular index pair. Indeed, we may construct a $C^{1}$  Lyapounov function $l:L_{1}\to
[0,1]$  for the flow $h^{f}$  with the additional property that  it vanishes precisely on
$L_{0}$.
Now, there is some $\tau>0$ such that  $l(x)>2\tau$ if $x\in S(P)\cup
S(Q)$. It now suffices to replace $L_{0}$ with $l^{-1}([0,\tau])$ and the pair
$(L_{1},L_{0})$ becomes a regular index pair. The fact that $I_{h^{f}}(K')=\emptyset$
implies that for each $x\in L_{1}-L_{0}$ the value
$t_{x}=sup\{t\in [0,\infty) : h^{f}_{[0,t]}(x)\subset L_{1}-L_{0} \}$
is well-defined and finite. Clearly, $t_{x}$ is the "arrival time" of $x$ in $L_{0}$.
Because the index pair $(L_{1},L_{0})$ is regular the arrival time function
$T:L_{1}\to [0,\infty)$ defined by $T(x)=t_{x}$ if $x\in L_{1}-L_{0}$ and
$T(x)=0$  otherwise, is continuous.  As $L_{1}$ is compact the function $T$
attains its maximum which, to simplify notation, can be assumed to be equal to $1$.

Let
$\psi :L_{1}\times [0,1]\to L_{1}$ be the function $\psi (x,t)=h^{f}_{min\{t,T(x)\}}(x)$.
This function is clearly continuous and, moreover, at a point  $(x,t)$ such that $t< T(x)$
it has the same order of differentiability as $h^{f}$. Let $\psi' :S(Q)\times [0,1]\to
L_{1}\times [0,1]$ be given by $\psi' (x,t)=(\psi(x,t),t)$. We denote by
$W$ the image of $\psi'$ and we let  $N'=Int(L_{1}-L_{0})\times [0,1]$ and $W'=W\cap
N'$. We  remark that $W'$  is a $C^{1}$ closed submanifold in $N'$. On the other hand
let $u:S(P)\to L_{1}$ be the inclusion. We consider the map $\phi:S(P)\times [0,1]\to
N'\times [0,1]$, $\phi (x,t)=(u(x),t)$. This map is $C^{1}$  and an embedding.
By standard transversality theory  we may find another $C^{1}$ embedding,
$\phi' :S(P)\times [0,1]\to L_{1}\times [0,1]$  which is arbitrarily close to $\phi$,
is transverse to $W'$ and coincides with $\phi$ when restricted to $S(P)\times\{0\}$. Let $S'$
be the image of $\phi'$ and notice that because $h^{f}_{T(x)}(x)\in L_{0}$ we obtain that
$S'\cap W'\cap (N'\times \{1\})=\emptyset$. If we let $C=(\phi')^{-1}(W')$ this implies that
$C$ is a $C^{1}$ submanifold of $S(P)\times [0,1]$ such that $\partial C=C\cap
(S(P)\times\{0\})=Z\times \{0\}$.

To see that $C$ does indeed provide the wanted null-bordism we still need to show
that the structure given by $(i,j)$ can be extended over $C$.  We start with the map
$j$ whose definition we recall. Let $w$ be an oriented path in $M$ that joins, in order, the
critical point $Q$ to the critical point $P$. For each point $y\in Z$ there is a unique flow
line, $w(y)$,  of the negative gradient flow $\gamma(f)$ of $f$ containing $y$ and joining
$P$ to $Q$.  The map $j$ is defined by $j(y)=w(y)\ast w$ where $\ast$ is concatenation of
paths (as $M$ is simply-connected the homotopy type of this map does not depend on $w$).
Notice that for each point $z\in S(P)$ there is a unique flow line $w'(z)$ of $\gamma(f)$
originating in $P$ and ending in $z$. Similarly, for each $z\in S(Q)$ there is a
unique flow line $w''(z)$ of $\gamma(f)$  originating in $z$ and ending in $Q$. If $y\in Z$,
then $w(y)=w'(y)\ast w''(y)$. Take now some point $(z,t)\in W'$. Then we have
$(z,t)= (h^{f}_{t}(x(z,t)),t)$ for a unique $x(z,t)\in S(Q)$.  Let
$w''(z,t)= \{x\in M :\exists v', \  f(Q)\leq v'\leq v\ , \ x=h^{f}_{t}(w''(x(z,t))\cap f^{-1}(v')\}$.
It is easy to see that
$w''(z,t)$ is a continuous path joining $(z,t)$ to $Q\times\{t\}$ and, of course,
$w''(z,0)=w''(z)\times\{0\}$.  Consider now a point $z'\in S'$. Then
$z'=\phi' (x'(z',t),t)$  for a unique couple $(x'(z',t),t)\in S(P)\times [0,1]$. We may assume
$\phi'$  close enough to the constant embedding $\phi$ such that for all $z'\in S'$ there is a
unique minimal geodesic in $V$ that joins  $(x'(z',t),t)$ to $z'$. We denote the path given by
this geodesic and ending in $z'$ by $w'''(z')$. Let $w'(z')=  (w'(x'(z',t))\times
\{t\})\ast w'''(z')$. This path starts in
$P\times\{t\}$ and ends in $z'\in S'$. If $(x,0)\in S(P)\times \{0\}$ then $w'(x,0))$ coincides
with $w'(x)\times \{0\}$.

We define $J':C\to  C^{0}(S^{1},M\times [0,1])$ by $J(z)=w'(\phi'
(z))\ast w''(\phi' (z))\ast w\times\{t\}$. It is very easy to check that this map is continuous.
The projection $M\times [0,1]\to M$ induces a coninuous map $e:C^{0}(S^{1},M\times
[0,1])\to C^{0}(S^{1},M)$ and therefore the composition $J=e\circ J'$ is continuous.
Notice also that the image of $J'$ is contained in $C^{0}((S^{1},0),(M,P))=\Omega(M)$.
It is clear that $J$ extends $j$.

We are now left to check that the framing associated to the embedding
$i:Z\hookrightarrow S(P)$ also extends to a framing of $C\subset S(P)\times [0,1]$.
We let $W(Q)=W^{s}(Q)\cap f^{-1}(-\infty, v]$. Of course $W(Q)$ is diffeomorphic
to a disk embedded as manifold with boundary in the pair $(f^{-1}(-\infty,v],V)$.
We fix the framing for the normal bundle of $W(Q)$. This means that we fix
an ordered family $(s_{1},  ..., s_{q})$ of linearly independent sections of this
bundle such that they provide the fixed  basis of $T_{Q}W(Q)$.
The framing on the normal bundle of $Z\subset S(P)$ asociated to $i$ is obtained
by projecting the sections $s_{i}$ on this normal bundle.  For $t\geq 0$ let
$W_{t}(Q)=h_{t}^{f}(W(Q))$. As $h^{f}_{t}$ is a diffeomorphism the
normal bundle to $W_{t}(Q)$ is well defined and we can transport by $(h^{f}_{t})_{\ast}$
the sections $s_{i}$ thus getting a framing of the normal bundle of $W_{t}(Q)$,
$(s_{1}^{t},\ldots,s_{q}^{t})$.  We now consider the manifold $\overline{W}\subset
M\times [0,1]$ defined by $\overline{W}=\cup_{0\leq t\leq 1}(W_{t}(Q),t)$. We may define
$q$ linearly independent sections in $TM\times [0,1]|_{\overline{W}}$ by
$(S_{i})_{(x,t)}=(s_{i}^{t})_{x}\times\{t\}$. We notice that this family of sections
$(S_{1},\ldots, S_{q})$ generate a supplement of $T\overline{W}$ in $T(M\times [0,1])$.
We now remark that $W'$ is a submanifold of $\overline{W}$.
Let $C'=W'\cap S'$. Of course, $C'=\phi' (C)$. Project the sections $S_{i}$ onto the
the normal bundle of $C'$ in $\phi' (S(P)\times [0,1])$. This produces $q$ linearly
independent sections of this bundle $(S'_{1}\ldots, S'_{q})$. We can transport these
sections back to $C$ and they give a framing of $C$ inside $S(P)\times [0,1]$ which
extends the framing of $Z\subset S(P)$ and cocludes the proof of the lemma.
\end{proof}

We now proceed to the proof of Proposition \ref{prop:bded_gen}

First we choose $\delta$ such that in the interval $(f(Q), f(Q)+\delta)$ there are no
critical values of $f|_{N}$  and $\delta\leq\delta'$  for $\delta'$ as in the proof of the last
lemma. We fix $v$ inside this interval and, as before we let
$V=f^{-1}(v)$. We now intend to adapt the proof of Lemma \ref{lem:bded_ho} to the present
situation.

We preserve all notations from the proof of that lemma.
The first remark  is that $S(Q)$ has the same properties as in that proof. However, we need to
replace $S(P)$ with  a different space that takes into account the fact that there might exist
critical points $R$ such that $f(Q)<f(R)<f(P)$. We let the space $\tilde{S}(P)$  consist of the
points  of $V$ that are situated on  possibly broken, flow lines that originate in $P$.  This space
is compact but it is not in general a manifold. It will be necessary to understand some of its
structure. For this let
$R_{1},\ldots, R_{h}$ be the critical points of $f$ such that $f(R_{i})\in (f(Q),f(P))$
and $R_{i}$ is connected to $P$ by a possibly broken flow line of $\gamma(f)$.
The Morse-Smale condition implies that $ind(R_{i})<p$. Therefore, the maximal
possible index of these critical points is $q$ (recall that by hypothesis $q$ and $p$
are succesive in $\{k\in {\bf N} : Crit^{N}_{k}(f)\not=\emptyset \}$). We will assume that
these critical points are ordered in decreasing order of their indexes so the first $s$ are all of
index $q$ and the rest are  of smaller indexes.  If no $R_{i}$ has index $q$ we take $s=0$.
The space $\tilde{S}(P)$ has the structure of a stratified space
$\tilde{S}(P)=\cup F^{k_{i}}$ with the strata $F^{k_{i}}$ of dimension $k_{i}-1$
defined as follows: $F^{k_{i}}$ is the space of the points $x\in V$ that are situated
on a non-broken flow line of $\gamma(f)$ originating in $P$ or in a critical point $R_{j}$
with $ind(R_{j})=k_{i}$. The first two top dimensional strata are therefore $F=F^{p}$ and
$F^{q}=\cup_{1\leq i\leq s}(W^{u}(R_{i})\cap V)$. We denote
$\overline{F}^{q}=\cup_{k_{i}\leq q}F^{k_{i}}$.

We now proceed as in the proof of Lemma \ref{lem:bded_ho}.
There is no modification in the construction of
the regular index pair $(L_{1},L_{0})$ only that we now ask that $\tilde{S}(P)$ be
contained inside $Int(L_{1}-L_{0})$.
The definition of $W'$ also remains the same. We let $\tilde{S'}$ to be the image of an
embedding $\tilde{\phi}'$ close to the constant embedding $\tilde{S}(P)\times [0,1]\subset
N_{1}\times [0,1]$ and such that $\tilde{\phi'}$ is transvere to $W'$ in the stratified sense
(this means that each stratum is itself transverse to $W'$). Of course, here
$\tilde{S}(P)\times [0,1]$ is stratified with the strata $F^{k_{i}}\times [0,1]$.
We denote $E^{k_{i}}=\tilde{\phi'}(F^{k_{i}}\times [0,1])$ and
$\overline{E}^{q}=\tilde{\phi'}(\overline{F}^{q})$.

If the intersection $W'\cap \overline{E}^{q}$ would be void, then
the intersection $C'=W'\cap \tilde{S'}$ would equal $D=W'\cap E^{p}$. Therefore,
it would  be a compact manifold providing a null-cobordism of $Z$ and the proof of
Lemma \ref{lem:bded_ho} would apply whithout modification.  However, $dim(W')=n-q$,
$dim(\overline{E}^{q})=q$ and transverality takes place inside $V\times [0,1]$ which is of
dimension $n$. This means that in general $W'$ might intersect
$\overline{E}^{q}$ non-trivially. At the same time, transversality implies that
$B=W'\cap \overline{E}^{q}=W'\cap E^{q}$ is a discrete union of points, $B=\{b_{1},\ldots
b_{r}\}$. Clearly, if we assume the condition (i) of the proposition, then
$s=0$, $B=\emptyset$ and this argument shows that the desired conclusion holds.

To show that this is also true under the assumption (ii) we need to analyze the general
situation. In this case $C'$ is a stratified manifold with (regular) boundary and with two
strata: the regular one, of dimension $p-q$,  and the second one equal  to  $B$.
Of course, $\partial{C'}=Z$ and the same method as that used in the proof of Lemma
\ref{lem:bded_ho} shows that we may extend the structure $(i,j)$ to $C'-B$.

Our purpose now is to show that each point of $B$ has a neighbourhood in $C'$ that is
a cone over its boundary. We will then delete from $C'$ the interior of a disjoint union of
neighbourhoods of this type, one for each point of $B$. This will then provide a bordism
between $Z$ and the disjoint union of the boundaries of these neighbourhoods.
But we will also see that each such boundary is null-bordant.

We will now focus on a single critical point $R_{i}$, $1\leq i\leq s$, such that
$ind(R_{i})=q$.  Let
$S^{u}(R_{i})=W^{u}(R_{i})\cap V$  and $S^{s}(R_{i})=W^{s}(R_{i})\cap
f^{-1}(f(R_{i})+\epsilon')$ such that
$\epsilon'$ is small enough such that there are no critical values in between $f(R_{i})$ and
$f(R_{i})+\epsilon'$. Of course, $S^{u}(R_{i})$ is a submanifold in $V$ but it is not compact
in general. Clearly, $S^{u}(R_{i})\subset F^{q}$. It is not difficult to see that
there is a neighbourhood $U(R_{i})$  of $S^{u}(R_{i})$ in $\tilde{S}(P)$ that is
homeomorphic to $S^{u}(R_{i})\times C(Z(P,R_{i}))$ where $C(-)$ is the un-reduced cone.
The key observation for this is that if one considers a tubular neighbourhood of
$S^{s}(R_{i})$ in $f^{-1}(f(R_{i})+\epsilon')$ then the boundary of this neighbourhood
intersects $W^{u}(P)$ along $S^{q-1}\times Z(P,R_{i})$. The gradient flow of $f$
carries a sphere of the form $S^{q-1}\times {pt}\subset S^{q-1}\times Z(P,R_{i})$ into
a copy of $S^{u}(R_{i})$ and therefore identifies the boundary of a neighbourhood of
$S^{u}(R_{i})$ in $\tilde{S}(P)$ to
$S^{u}(R_{i})\times Z(P,R_{i})$. Let now $U'(R_{i})=\tilde{\phi'} (U(R_{i})\times [0,1])$.
Clearly, $U'(R_{i})$ is homeomorphic to $(S^{u}(R_{i})\times [0,1])\times
CZ(P,R_{i})$.  We also have $E^{q}=\cup_{1\leq i\leq s} \tilde{\phi}'(S^{u}(R_{i})\times
[0,1])$. Assume that the point $b_{j}\in W'\cap E^{q}$ belongs to
$\tilde{\phi'}(S^{u}(R_{i}))$. This means that $b_{j}\in W' \cap
\tilde{\phi}'(S^{u}(R_{i})\times [0,1]))$ and  this last intersection is transverse. By
trnasversality, the intersection $U'(R_{i})\cap W'$ contains a neigbourhood $U_{b_{j}}$ of
$b_{j}$ in $W'$ which is homeomorphic to $CZ(P,R_{i})$.

By eliminating from $C'$ the disjoint union of the interiors of $U_{b_{j}}$
we obtain a manifold $C''$ with boundary such that $\partial C''$ is the disjoint
union of $Z$ and a number of copies $Z^{j}_{i}=\partial U_{b_{j}}$ of some $Z(P,R_{i})$,
$1\leq i\leq s$. As discussed above the structure $(i,j)$ extends to $C''$.
The purpose now is to show that the restriction of this structure
to $\coprod Z^{k}_{i}$  produces a bordism class which is the same
with the sum of the bordism classes given by $j(P,R_{i}):Z(P,R_{i})\to \Omega M$
together with the framings associated to $Z(P,R_{i})\subset S^{u}(P)$.  As all the
stabilizations $[H(j_{P}, j_{R_{i}})]$ of the Hopf invariants  are vanishing by hypothesis,
this implies that
$Z$ is null-bordant  and therefore, by Theorem \ref{theo:hopf_conn} again,
$[H(j_{P},j_{Q})]$  is null thus contradicting the hypothesis.

Let $J$ be the extension
of the map $j$ to $C''$ constructed as in the proof of Lemma \ref{lem:bded_ho}. It is
easy to notice that $J|_{Z^{k}_{i}}\he j(P,R_{i})$ because, with the notations
in the proof or \ref{lem:bded_ho},  we have that for all $z\in Z^{j}_{i}$,
$w''(z)$ is arbitrarily close to a constant path joining $b_{j}$ to $Q$. On the other hand
$w'(z)$ is very close a path joining $P$ to $R_{i}$ and passing through $z'$ where $z'$ is
identified to $z$ via the homeomorphism $Z^{j}_{i}\approx Z(P,R_{i})$  and  followed
by a flow line from $R_{i}$ to $b_{j}$.

The next task is to discuss the framings (this is the only place the non-vanishing of
$[H(j_{P},j_{Q})]$ is needed and that of  $H(j_{P},j_{Q})$ is not sufficient).
For this consider
$U(\overline{E}^{q})\subset
\tilde{S}'$ to be a closed neighbourhood of the singular set $\overline{E}^{q}$ of
$\tilde{S}'$ of the from $\tilde{\phi}'(U\times [0,1])$ with $U$ a neighbourhood of the
singular set of $\tilde{S}(P)$ and such that $W'\cap (\tilde{S}'-U(\overline{E}^{q}))=C''$,
$\coprod Z^{j}_{i}=W'\cap \partial U(\overline{E}^{q})$ and $\coprod U_{b_{j}}=W'\cap
U(\overline{E}^{q})$. Write the intersection point $b_{j}\in W'\cap \tilde{S}'$
as $b_{j}=(a_{j},t_{j})\in V\times [0,1]$. By slighly modifying $W'$ (using a compactly
supported isotopy) in the neighbourhood of these intersection points we may assume that the
$t_{j}$'s are pairwise distinct  and that  $U_{b_{j}}\subset V\times \{t_{j}\}$.
We assume that the $t_{j}>t_{j-1}$ for all the $b_{j}$'s.
Let $\delta>0$ be smaller than $t_{1}$, $1-\delta > t_{r}$ and consider the manifold defined
as $S'' =(\tilde{S}'-Int(U(\overline{E}^{q})))\cup_{\partial U(\overline{E}^{q})} D\times
[0,\delta]$ where $D\subset S^{p-1}$ is such that $(\tilde{S}(P)-U)\cup D=S^{p-1}$.
The existence of $D$ is verified by using the inverse flow of $f$ to move $\tilde{S}(P)-U$
inside $W^{u}(P)\cap f^{-1}(P-\epsilon'')$ for $\epsilon''$ very small. We can take
$D$ to be the closure of the complement of the image of $\tilde{S}(P)-U$ inside this
$p-1$-dimensional sphere. We now remark that, after rounding the corners of $S''$,
we have that $S''\approx S^{p-1}\times [0,1]$ via a homeomorphism such that
$Z\subset \partial_{1}S''\approx S^{p-1}\times \{0\}$ and $\coprod Z^{j}_{i}\subset
\partial_{2}S''\approx S^{p-1}\times \{1\}$.
As in the proof of Lemma \ref{lem:bded_ho} we can construct $q$ linearly independent
sections $(S'_{1},\ldots, S'_{q})$ of a supplement of $TW'$ in $T(V\times [0,1])$.
These sections provide the framing of $Z$ in $\partial_{1}S''$. They provide
a framing of $\coprod Z^{j}_{i}$ in $\partial_{2}S''$ and make $C''$
a framed cobordism inside $S''$. The key fact to recall now is that
in general, unstably, the framing of a disjoint union is not equivalent to the sum
of the framings of the components (two framings are equivalent if they are related by
a framing of the product of the manifold with the unit interval). However, stably this is
true - the (stable) bordism class of any disjoint union is just the sum of the bordism classes of
the components. Therefore, we can now study the framing of each
$Z^{j}_{i}$ at a time and, to conclude the proof of the proposition, it will be enough to show
that the framing $\xi$ induced on $Z^{j}_{i}$ by the sections $(S'_{k})$ is equivalent at least
stably with the standard framing of $Z(P,R_{i})\subset S^{u}(P)$  (stable equivalence
meaning that the two  framings are equivalent after possibly adding a trivial bundle with a
fixed framing to the two normal bundles in question).  To simplify notation
we may assume that in a neighbourhood
$\mathcal{N}$ of $b_{j}$ in $V\times [0,1]$ we have
$\tilde{\phi'}(S^{u}(R_{i})\times [0,1])\cap\mathcal{N}=S^{u}(R_{i})\times
[0,1]\cap\mathcal{N}$. So for this argument we may assume
$\tilde{\phi'}(S^{u}(R_{i})\times [0,1])=S^{u}(R_{i})\times [0,1]$.
Recall that we have the inclusion
$Z^{j}_{i}\hookrightarrow CZ^{j}_{i}=U_{b_{j}}$. Clearly, there is an obvious map
$CZ^{j}_{i}\to b_{j}$ and, because the intersection of $W'$ with
$S^{u}(R_{i})\times [0,1]$  is  transverse in $b_{j}$ we have that the sections
$(S'_{k})$ project to a basis in $T_{b_{j}}(S^{u}(R_{i})\times [0,1])$  which
can be viewed as a section of the normal bundle of $b_{j}$ in $S^{u}(R_{i})\times [0,1]$ .
By looking to the composition $Z^{j}_{i}\hookrightarrow CZ^{j}_{i}\to b_{j}$ we see that
the framing  of $Z^{j}_{i}$  is induced from this basis by this map (in the
sense that we have maps between the total spaces of the normal bundles of the three
spaces considered that make the sections correspond one to the other).
Moreover,  the contractibility of $CZ^{i}_{j}$ implies that any two framings
induced in this way from two possibly different basis are equivalent.

We now need to use the relation between $Z^{j}_{i}$ and $Z(P,R_{i})$.
We fix $q$ sections $(x_{1},\ldots,x_{q})$ of the tangent bundle of $W^{u}(R_{i})$.
We can also extend them in the obvious way to a neighbourhood of $R_{i}$
contained in $W^{u}(R_{i})\times W^{s}(R_{i})$.
The value of these sections in $R_{i}$ induces the normal framing of $W^{s}(R_{i})$
and that of $Z(P,R_{i})$ in $S^{u}(P)$.

We recall
that  $\{a_{j}\}\times Z^{j}_{i} \subset S^{u}(R_{i})\times
(CZ^{j}_{i}-{\ast})\subset (S^{u}(R_{i})\times D^{n-q})\cap W^{u}(P)$
is obtained (by flowing via the negative gradient of $f$) from
$\{a_{j}\}\times Z(P,R_{i})\subset  S^{q-1}\times Z(P,R_{i})=\partial D^{q}\times
Z(P,R_{i}\subset W^{u}(P)\cap S^{s}(R_{i})$ (of course all this operation takes place
in $M\times\{t_{j}\}$ but we will neglect this parameter).  Fix  the framing $\xi'$
of the normal bundle of  $\{a_{j}\}\times CZ^{j}_{i}$
inside $\overline{(W^{u}(R_{i})\times D^{n-q})\cap W^{u}(P)}$
induced by the sections $\{(x_{k})\}$. It only depends on the value of these sections in the
point $\{a_{j}\}\times \{\ast\}\in S^{u}(R_{i})$. At the same time
the relation between $Z^{j}_{i}$ and $Z(P,R_{i})$  means
that $\xi'$ is equivalent to the standard framing of $Z^{i}_{j}$ in $S^{u}(R_{i})$.

We are now left to compare the framings $\xi$ and $\xi'$ of the normal bundles
of $CZ^{j}_{i}$ associated, the first to the stratified embedding
$e_{1}:(CZ^{j}_{i},\ast)\to (\partial_{2}S''=S^{p-1}\times [0,1], S^{u}(R_{i})\times [0,1])$
and the second to the embedding $e_{2}:(CZ^{j}_{i},\ast)\to
(\overline{W^{u}(P)},W^{u}(R_{i}))$. Both these two embeddings factor via the same
embedding $(CZ^{j}_{i},\ast)\to (\tilde{S}(P),S^{u}(R_{i}))$ and the framings are induced
in the way described before by pull-back over the composition $Z^{j}_{i}\hookrightarrow
(CZ^{j}_{i},\ast)\to\ast$. After adding a rank one trivial line bundle to the two normal bundles
in question we see that $\xi$ and $\xi'$ become equivalent.
This concludes the proof of the proposition.
\end{proof}

\

{\em Proof of Theorem \ref{theo:closed}.}

\

We assume that $f$, $P$ and $Q$ are as in the statement. As $[H(j_{P},j_{Q})]\not=0$
there is some critical point $Q'\in Crit^{N}_{q}(f)$ such that $[H(j_{P},j_{Q'})]\not=0$,
$f(Q')\geq f(Q)$ and for all $R\in Crit^{N}_{q}(f)$, $f(Q')<f(R)<f(P)$ we have
$[H(j_{P},j_{R})]=0$. We can therefore apply
Proposition \ref{prop:bded_gen}  to the pair $P,Q'$. We obtain that there exists
some $\delta>0$ such that there are infinitely many regular values
$v_{i}\in(f(Q)+\delta,f(Q)+2\delta)$  (of $f|_{N}$) such that each hypersurface
$V_{v_{i}}=f^{-1}(v_{i})\cap N_{1}$ contains at least one bounded orbit of $h^{f}$.  We
may assume that the set $\{v_{i}\}$ is discrete. We fix one such $v_{i}$. Consider a point
$x\in V_{v_{i}}$ whose orbit is bounded. Then $\omega (x)$ is non-void. Let $y\in\omega
(x)$.Then $y$ is non-wandering and belongs to $\Omega_{c}(h^{f})$. We now fix some
small $\epsilon_{i}>0$ and apply the $C^{1}$ closing  lemma of Pugh and Robinson
\cite{PuRo} to deform $H_{f}$ by a small $C^{1}$ deformation with a fixed
compact support (only depending on the orbit to be closed) contained in
$f^{-1}(v_{i}-\epsilon_{i},v_{i}+\epsilon_{i})$ to obtain a hamiltonian vector field
$H_{f'}$ such that $h^{f'}$ has a periodic orbit close to $y$. Because the deformation used has
a fixed compact support, $f'$ can be made as $C^{2}$-close to $f$ as needed and such that
the constructed periodic orbit is included in
$f^{-1}(v_{i}-\epsilon_{i},v_{i}+\epsilon_{i})$. We can then apply the same process
for all the values $v_{i}$ by taking the $\epsilon_{i}$'s such that the intervals
$f^{-1}(v_{j}-\epsilon_{j},v_{j}+\epsilon_{j})$ are pairwise disjoint,
$\lim_{j\to\infty}\epsilon_{j}=0$ and, when
$j\to\infty$, the deformations used tend to $0$. This leads to the  existence of the
function $\overline{f}$ and concludes the proof.

\subsection{Applications and examples.}\label{subsec:appl}

\

The first paragraph below is concerned with a
disscussion of property $(\ast)$. The main part of the sub-section appears in the second
paragraph which  contains corollaries of  Theorem \ref{theo:closed}. We end with some open
problems and further questions.

\subsubsection{Property $(\ast)$.}

\

A few simple critera that insure the validity of property $(\ast)$ are contained in the following
lemma.

\begin{lemma}\label{lem:prop_ast}
Let $f:M\to {\bf R}$ be a $C^{2}$ function. Let $S$ be an isolated invariant set for the
negative gradient flow $\gamma=\gamma(f)$ of $f$.
The invariant set $S$ satisfies property  $(\ast)$ if one of the following conditions
is satisfied.
\begin{itemize}
\item[(i)] There exist an isolating neighbourhood $N$ of $S$, $v\in{\bf R}$ regular for
$f|_{Int(N)}$  and a neighbourhood  $U\subset N$ of $W^{s}_{N}(S)$  such that
$v< inf\{ f(x):x\in S \}$ and, for all $x\in U-W^{s}_{N}(S)$, there is $t\in {\bf R}^{+}$ with
$f(\gamma_{t}(x))=v$, $\gamma_{[0,t]}(x)\subset N$.
\item[(ii)] There are regular values $a<b$ of $f$ such that $S$ is the maximal compact invariant
set contained in $f^{-1}([a,b])$ and for any $x\in f^{-1}([a,b])$ we have: $\omega(x)\subset S$
or there exists some $t\in{\bf R}^{+}$ such that $f(\gamma_{t}(x))=a$ and
$\omega^{\ast}(x)\subset S$ or there is $t\in {\bf R}^{+}$, $f(\gamma_{-t}(x))=b$.
\item[(iii)] $S$ is the maximal compact invariant set of $\gamma$, there exists a compact $K\in
M$ and an $\epsilon>0$ with $||\nabla f(x)||>\epsilon$ for $x\not\in K$ and $\gamma$ is a flow
(not only a partially defined one).
\item[(iv)] $S$ is the maximal compact invariant set of $\gamma$, $M$ is metrically complete
and there exist $m>\epsilon>0$ such that outside some compact set we have
$m>||\nabla f(x)||>\epsilon$.
\end{itemize}
\end{lemma}
\begin{proof}
It is clear how to use (i) to prove property $(\ast)$: we consider $U'\subset Int(U)$
another neighbourhood of $W^{s}_{N}(S)$ and we take $N_{1}=\{x\in N: \exists t >0 \ , \
\gamma_{[-t,0]}(x)\subset N \ ,\ \gamma_{-t}(x)\in \overline{U'}\}$ and we let
$N_{0}=N_{1}\cap f^{-1}(v)$.
From (ii) we immediately deduce that we may construct $N$, $U$ as required at (i).
Point (iii) implies that the condition at (ii) is verified; (iv) implies (iii).
\end{proof}

\begin{remark}
It is immediate that condition (i) in the lemma above is, in fact equivalent to property $(\ast)$.
\end{remark}

It is easy to produce examples when condition $(\ast)$ is not satisfied.

\begin{example} \label{exam:non_prop} Let $g:D^{2}\to {\bf R}$ be a Morse-Smale function
with  precisely two critical points: $0$ which is of index $2$ and $Q\in Int(D^{2})$ of index
$1$. We may also find $g$ such that  $g(0)=1$, $g(Q)=0$ and there is precisely one (negative
gradient) flow line $\xi$ in $D^{2}$ that joins $0$ to $Q$. We consider
$D'=D^{2}-\{T\}$ where $T$ is a point in $D^{2}$ which is not in $\xi\cup \{0,Q\}$ and
$g(T)>0$. We let $f=g|_{D'}$ and take $S=\{0,Q\}\cup \xi$.  It is easy to check that $S$
does not have the  property $(\ast)$.
\end{example}

A simple result concerning this property $(\ast)$ will be of use later.
We recall that $\x(M)$ is the space of $C^{2}$ local isolated functions on $M$,
that is pairs $(f,N)$ with $f$ a $C^{2}$ function and $N$ an isolating neighbourhood
of the negative gradient flow $\gamma(f)$ of $f$.

\begin{lemma}\label{lem:stab_ast} Let $(f,N)\in\x(M)$ and $S=I_{N}(f)$. Suppose
that $S$ verifies property $(\ast)$. There is a $C^{2}$-neighbourhood $\mathcal{U}$ of $f$
such that if $g\in\mathcal{U}$, then $(g,N)\in\x(M)$ and $I_{N}(g)$ verifies $(\ast)$.
\end{lemma}

\begin{proof}
Fix $(N_{1},N_{0},N_{2})$ a strong index block of $S$ with $N_{0}\subset f^{-1}(a)$
and $f$ regular on $N_{0}$. It is clear that by taking $\mathcal{U}'$ a sufficiently small
neighbourhood of $g$ we have that $(N_{1},N_{0},N_{2})$ is also
a strong index block for each $g\in\mathcal{U}'$.  Fix a regular value of $f|_{N_{1}}$,
$a'>a$ such that $a'<inf\{f(S)\}$. It is easy to see that by possibly diminishing more
$\mathcal{U'}$ the pair $(N_{1}\cap g^{-1}([a',\infty)),N_{1}\cap g^{-1}(a'))$
becomes an index pair for each such $g$.
\end{proof}

\subsubsection{Corollaries of  Theorem \ref{theo:closed}.}

\

We start by defining the invariant $d(-)$ mentioned in the introduction.
This invariant is inspired by the fact that
in the statement of  Theorem \ref{theo:closed} the Hopf invariant appearing
there is a  coefficient in the differential of the extended Morse complex
$\mathcal{C}_{N}(f)$.

\begin{definition}\label{def:inv} Let $(f,N)\in\x(M)$.
 If $(f,N)$ is a local Morse function
 we  let $d_{N}(f)\in {\bf Z}/2$ be equal to $1$ iff $\mathcal{C}_{N}(f)$ has a non-trivial
differential.  For a general $(f,N)\in\x(M)$ we let $d_{N}(f)=sup\{inf\{d_{N}(g)
:g\in\mathcal{U},\  (g,N)
\ local\ Morse \}:
\mathcal{U} \  open \ in \ C^{2}(M,{\bf R}), f\in \mathcal{U}\}$ (here $C^{2}(M,{\bf R})$ is
the space of $C^{2}$ functions with the strong topology).
\end{definition}

Clearly, if $d_{N}(f)=1$, then all Morse functions $g$
which are sufficientely close to $f$ have an extended Morse complex $\mathcal{C}_{N}(g)$
with non-trivial differential. Moreover, there is a neighbourhood of $f$
consisting of functions $f'$ that have $N$ as an isolating neighbourhood and with
$d_{N}(f')=1$.

The local isolated functions $(f,D)$ with $f$
appearing in Example \ref{exam:non_bded} and $D\subset {\bf R}^{2k}$ a closed disk
containing the origin verify $d_{D}(f)=0$.

Assume now that $f$ is a $C^{2}$ function and $S$ is an isolated invariant set
for $\gamma(f)$. We let $d(S)=d_{N}(f)$ where $N$ is an isolating neighbourhood
such that $S=I_{N}(f)$. It is immediate that $d(S)$ is indeed independent of $N$.
We also let $e(S)$ to be the infimum of  $d_{N'}(f')$ where the pair $(f',N')\in\x(M)$ is
continuation equivalent to $(f,N)$ for some $N$ with $S=I_{N}(f)$. Assuming
$(f,N)\in\x(M)$ fixed, notice, as above, that a function $f'$ that is sufficiently close to $f$
has $N$ as an isolating neighbourhood and $(f',N)$ is continuation equivalent to $(f,N)$.
This means that $d(S)\geq e(S)$. It is very easy  to construct examples
such that $d(S)>e(S)$. For example, take $S=\{0,Q\}\cup \xi$ in
Example \ref{exam:non_prop}.

\

Here is an important consequence of Theorem \ref{theo:closed}.

\begin{corollary}\label{cor:inv} Let $(f,N)\in\x(M)$ such that $S=I_{N}(f)$ satisfies property
$(\ast)$ and $d(S)=1$.
\begin{itemize}
\item[(i)] If $(f,N)$ is a local Morse function, then there are infinitely many regular
hyersurfaces $V=f^{-1}(v)$ which, after perturbation by an arbitrarily small compactly
supported isotopy of $M$ carry a closed characteristic.
\item[(ii)] There is a family of functions, $\mathcal{F}$,  dense
in a $C^{2}$ neighbourhood of $f$ such that each hamiltonian flow $h^{f'}$,
$f'\in\mathcal{F}$  has infinitely many distinct periodic orbits.
\end{itemize}
\end{corollary}

\begin{proof}
We prove point (i). We first show the existence of one hypersurface $V$.
If $(f,N)$ is a local Morse function and $\mathcal{C}_{N}(f)$ has a non-trivial
diferential we obtain that there is a pair of  consecutive critical points $P,Q\in Crit^{N}(f)$
verifying the assumption of Theorem \ref{theo:closed}. By the proof of this theorem
we obtain that  we may find a regular value $v$ of $f$ such that in
any $C^{2}$ neighbourhood of $f$ and for any $\epsilon$ there is a function $\overline{f}$
whose associated hamiltonian has inside $Int(N)\cap \overline{f}^{-1}(v-\epsilon,v+\epsilon)$
at least one periodic orbit.  Assume that this periodic orbit actually belongs to
$\overline{f}^{-1}(v')$. By taking $\epsilon$ small enough and $\overline{f}$ very
close to $f$ it is easy to see that there is a compactly supported isotopy  $\phi$ of $M$ such
that $\phi(V)\cap Int(N)=\overline{f}^{-1}(v')\cup Int(N)$. This means that the hypersurface
$V'=\phi(V)$ carries a closed characteristic. Now the same argument can be applied for
infinitely many hypersurfaces $V$ as in the proof of Theorem \ref{theo:closed}.

The second point is essentially obvious given
Theorem \ref{theo:closed}.  By Lemma \ref{lem:stab_ast}
all the functions in a sufficiently small neighbourhood $\mathcal{U}$ of $f$
are such that $g\in \mathcal{U}$ implies $(g,N)\in\x(M)$ and $I_{N}(g)$
verifies property $(\ast)$.  As seen above me may assume $\mathcal{U}$
small enough such that each $g\in\mathcal{U}$ also satisfies $d(I_{N}(g))=1$.
Local Morse functions form a dense family in $\mathcal{U}$. For such a function
$(h,N)$ we have $d_{N}(h)=1$ which means that $\mathcal{C}_{N}(h)$ has a non trivial
differential. But by Theorem \ref{theo:hopf_conn} this means that the assumptions
needed to apply Theorem \ref{theo:closed} are satisfied. Therefore arbitrarily close to
$h$ there is a function $h'$ whose hamiltonian flow has infinitely many periodic orbits.
This proves the claim at (ii).
\end{proof}

\begin{remark} It is useful to note that the  hypersurfaces produced
at the point (i) of  Corollary \ref{cor:inv} might all be isotopic one to the other.
\end{remark}

It is hard in general to determine $d(f)$. Of course, the exact determination of $e(f)$
is not simple either. However, and this is fundamental for applications, we have
a reasonably effective criterion that insures $e(f)=1$.

\begin{corollary} \label{cor:cond} Let $(f,N)\in \x(M)$ and $S=I_{N}(f)$. Each one of the
following conditions implies that $e(S)=1$.
\begin{itemize}
\item[(i)] There exists $(g,N')$
 which is local Morse, continuation equivalent to $(f,N)$  and such that the differential
$d_{g}$ in $\mathcal{C}(g)$ verifies $d_{g}\not=d^{Mo}\otimes id_{\Omega M}$.
\item[(ii)] There is no minimal local Morse function that is continuation equivalent
to $(f,N)$
\item[(iii)] There is a minimal Morse function $(g,N')$ continuation equivalent to $(f,N)$
and its extended Morse complex
$\mathcal{C}_{N'}(g)$ has a non-trivial differential.
\item[(iv)] $H_{\ast}(c_{N}(f);{\bf Z})$ has torsion.
\item[(v)] $H_{\ast}(c_{N}(f);{\bf Z})$ is torsion
free, $\overline{c}_{N}(f)$ is simply-connected and there are $k,l\in{\bf
N}$ with $H_{t}(c_{N}(f);{\bf Z})=0$ for $0<k<t<l$, $l-k>1$ and in the homology Serre
(or in the $\pi_{\ast}^{S}(\Omega M)^{+})$- Atiyah-Hirzebruch-Serre) spectral sequence of
the fibration induced by the canonical projection (see \S\ref{subsubsec:conlind})
$\overline{c}_{N}(f)\stackrel{p}{\to} M$  from
$\Omega M\to PM\to M$  the differential $d^{l-k}$ (resp. $D^{l-k}$) does
not vanish when restricted to
$E^{l-k}_{l\ast}$  (resp. $\mathcal{E}^{l-k}_{l\ast}$).
\end{itemize}
\end{corollary}

\begin{proof}
For any local Morse function $(g,N')$ the differential $d^{g}$ of the complex
$\mathcal{C}_{N'}(g)$ contains $d^{Mo}$  in the sense that $(d^{Mo}\otimes id_{\Omega
M})(x)=d^{g}(x)$ for all critical points $x$ such that $Crit^{N'}_{|x|-1}(g)\not=\emptyset$.
Therefore, $d^{Mo}\not=0$ implies $d^{g}\not=0$. This means that if $g$ is not minimal,
then $d(g)=1$. Therefore, to decide that $e(S)=1$ we only need to make sure that each minimal
local Morse function
$(g,N')$  which is continuation equivalent to $(f,N)$ verifies $d_{N'}(g)=1$. But by Corollary
\ref{cor:PS_min} any two minimal Morse functions that are continuation equivalent have
isomorphic extended Morse complexes.  Thus it is enough to verify this condition for just one
minimal local Morse function  that is continuation equivalent to $(f,N)$.

This remark directly justifies our claim if one of the conditions (ii),(iii) or (iv) is satisfied.
It also implies the claim if (v) is satisfied. Indeed, asuming (v) and supposing that a
minimal local Morse function $(g,N')$  is continuation equivalent to $(f,N)$
 its Morse complex $\mathcal{C}_{N'}(g)$ is isomorphic by Corrolary \ref{cor:min} to the
chain complex  $\mathcal{C}_{min}(\overline{c}_{N}(f) )$ (which is defined under the
assumptions made). At the same time  Corollary \ref{cor:s_ss} shows that the differential of
this last complex is non-trivial.

We are now left to prove that the point (i) also implies $e(S)=1$.
The condition at (i) implies that
there is some critical point $P\in Crit^{N'}(g)$ whose differential in the classical Morse
complex is null but its differential in $\mathcal{C}_{N'}(g)$ is not vanishing. But this means
that $P$ represents a non-trivial homology class in the homology of the classical Morse
complex and that there is another critical point $Q\in Crit^{N'}(g)$  such that
$[H(P,Q)]\not=0$. It also results that if $q=ind_{g}(Q)$ and $p=ind_{g}(P)$, then $p>q+1$.
As explained before we need to consider the case when there exists a minimal
local Morse function $(h,N'')$ continuation equivalent to $(f,N)$ (and thus to $(g,N')$)
and it is enough to show that $d_{N''}(h)=1$.  We obtain by Corollary \ref{theo:rig} that there
exists a map $m:C^{Mo}_{N'}(g)\to C^{Mo}_{N''}(h)$
which induces a morphism of extended Morse complexes for all compatible index sets.
Let $I=\{0,1,\ldots, q-1,q,p,p+1,\ldots , n\}$. Then both complexes
$\mathcal{C}^{I}_{N'}(g)$ and $\mathcal{C}_{N''}^{I}(h)$ are defined. Denote by $m^{I}$
the morphism induced by $m$ between these two extended complexes. We denote
by $\tilde{d}$ the differential in these two complexes.
As $m$ induces an isomorphism in homology and $P$ represents a non-trivial
homology class in $H_{\ast}(C^{Mo}(g))$ it results that there is $u\in
{\bf Z}[Crit^{N''}_{p}(h)]$ such that $m(u)=P+d^{Mo}(b)$ for some $b\in
C^{Mo}_{N'}(g)$ and $d^{Mo}(u)=0$ because $p>q+1$ and $h$ is minimal.  Now, we have
that $\tilde{d}d^{Mo}=0$. Therefore $\tilde{d}(m^{I}(u))=[H(j_{p},j_{Q}]Q+c$ where $c$
does not contain any terms in $Q$. Therefore, $\tilde{d}(m^{I}(u))\not=0$. But this implies
that $\tilde{d}u\not=0$. Therefore the differential $\tilde{d}$ is not zero in
$\mathcal{C}^{I}_{N''}(h)$. But the  differential in $\mathcal{C}_{N''}(h)$ contains
$\tilde{d}$ and therefore is also non-trivial.
 \end{proof}

\begin{remark}\label{rem:var}
a. The conditions (iv), (v) appearing in Corollary \ref{cor:cond} are continuation invariant.
However, sometimes it is easier to apply a version of these
which uses instead of the pair $\overline{c}_{N}(f)$ an index pair $(N_{1},N_{0})$
of the invariant set $I_{N}(f)$.
In this case condition (iv) remains the same for recall that $H_{\ast}(c_{N}(f);{\bf
Z})=H_{\ast}(N_{1},N_{0};{\bf Z})$.  At point (v) the
assumption that replaces the simple-connectivity of $\overline{c}_{N}(f)$ is that
both $N_{1}$ and $N_{0}$ are simply connected. Of course, instead
of the spectral sequences associated to the pair $\overline{c}_{N}(f)$ one may use
the same spectral sequences for the pair $(N_{1},N_{0})$
Finally, in view of Corollary \ref{cor:s_ss}, this spectral sequence condition as it appears at
point (v) is equivalent with the fact that the minimal Hopf complex of the pair
$(N_{1},N_{0})$ has a non trivial differential.

b. It is easy to find manifolds $M$ and Morse functions $f:M\to {\bf R}$ such that condition (i)
of the Corollary \ref{cor:cond} is satisfied (of course, $\pi_{\ast}^{S}(\Omega M^{+})$ needs
to be nontrivial for $\ast>0$). A simple example is provided by the height function $h$ on
a sphere $S^{k}$ for $k\geq 1$. This has precisely two critical points: $T$  a
maximum and $S$ a minimium. The differential in the classical Morse complex is null.
However, in $\mathcal{C}_{S^{k}}(h)$ we have $dT=(Id)^{\ast}$ where $(Id)^{\ast}$ is the
stable homotopy class of the inclusion of the bottom cell $S^{k-1}\hookrightarrow \Omega
S^{k}$. Another instructive example is a minimal Morse  function $\tilde{h}$ on ${\bf
C}P^{n}$.  Again, the classical Morse differential vanishes but on each critical point
 $P$ of $\tilde{h}$ with $ind(P)>0$  the differential in the extended Morse complex
$\mathcal{C}(\tilde{h})$ is not null. Of course, both these examples are compact but
via a product with some non-compact manifold (${\bf R}^{k}$, for example) or,
 more generally, by considering some
vector bundles over them, they  produce non-compact variants.
\end{remark}

The Serre spectral sequence condition at (v) is in many cases easy to verify.

\begin{example}\label{ex:cp} Let $M={\bf C}P^{n}\times {\bf R}^{2k}$ and
consider the function $g\oplus q$ where $g:{\bf C}P^{n}\to {\bf R}$ is a perfect Morse function
and $q:{\bf R}^{2k}\to {\bf R}$ is a quadratic form of index different from
$0$ or $2k$. Fix on
$M$ an arbitrary symplectic form and let $f$ be obtained by a compactly supported
perturbation of $g\oplus q$
(notice that if we take on $M$ a symplectic form that is simply the sum of the
forms on ${\bf C}P^{n}$ and ${\bf R}^{2k}$, then it is obvious that $h^{g\oplus q}$ has
many bounded orbits because the hamiltonian flow leaves invariant
${\bf C}P^{n}\times \{0\}\subset M$).

The gradient flow of $f$ has a maximal compact invariant set $I(f)$
and it is easy to see that this set satisfies the condition $(\ast)$ (this follows by applying
(ii) of Lemma \ref{lem:prop_ast}). Notice also that $d(I(f))=1$. Indeed, there exists an index
pair of $I(f)$ of the form $(N_{1},N_{0})=
({\bf C}P^{n}\times D^{2k}, {\bf C}P^{n}\times S^{j-1}\times D^{2k-j})$
where $j=ind(q)$. Consider the homology Serre spectral sequence $(E^{r}_{pq},d^{r})$ of the
fibration $\Omega {\bf C}P^{n}\to (E_{1},E_{0})\to (N_{1},N_{0})$ which is obtained by
pull-back of the (homotopy) fibration $\Omega {\bf C}P^{n}\to \ast \to {\bf C}P^{n}$.
It is immediate to see that this spectral sequence is simply a suspension of the
Serre spectral sequence $(e^{r}_{pq},\delta^{r})$ of this last fibration in the
sense that $E^{r}_{(p+j)q}=e^{r}_{pq}$ and similarly for the differentials.
In turn this fibration has all differentials $\delta^{2}_{(2l)0}$ non trivial. This means
that the conditions needed to apply Corollary  \ref{cor:cond} (v) are satisfied.

This example also indicates that the non-vanishing of $d(-)$ is sometimes implied
by a non-trivial cup-structure in co-homology (this is of course to be expected
as the non-vanishing of Hopf invariants is implied by the existence
of certain non-trivial cup-products).
\end{example}

Some of the conditions in the last corollary imply slightly more than what is
claimed there. The next result also emphasizes the stability properties of the extended
Morse complexes.

\begin{corollary} \label{cor:cl_ht} Let $(f,N)\in \x(M)$ and let $(g,N')\in\x(M\times {\bf
R}^{2k})$ such that $(g,N')$ is continuation equivalent
to $(f\oplus q,N\times D^{2r})$  with $q:{\bf R}^{2r}\to {\bf R}$ a
non-degenerate quadratic form.
If one of the conditions (iv) and
(v) from Corollary
\ref{cor:cond} are satisfied for $(f,N), S=I_{N}(f)$, then there exists a family of
functions  dense in a neighbourhood of $g$  each of whose induced hamiltonian
flow on $M\times {\bf R}^{2r}$  (endowed with an arbitrary symplectic form) has
infinitely many distinct  periodic orbits.
\end{corollary}
\begin{proof}
Condition (iv) shows that there are no minimal
Morse functions in the continuation equivalence class of the sum $(f\oplus q,N\times
D^{2r})$ (because the Conley index of the maximal invariant set associated to such a function
is just a suspension of the Conley index asssociated to $(f,N)$ ).
Similarly,  by the definition of the global Conley index, we see that
$\overline{c}_{N\times D^{2r}}(f\oplus q)$ is simply connected and torsion free and therefore
the minimal Hopf complex  $\mathcal{C}_{min}(\overline{c}_{N\times D^{2r}}(f\oplus q))$
is defined. It is  easy to see that this complex is just the suspension of order $ind(q)$ of
$\mathcal{C}_{min}(\overline{c}_{N}(f))$. Thus, as condition (v) implies
that the differential in this last complex is not trivial, we obtain that the differential of the first
is not trivial either. This immediately implies the claim by using Corollary  \ref{cor:cond}.
\end{proof}

We now present a class of  examples where the methods decribed above
apply naturally.

\begin{example}\label{exam:appli}
 Let  $B$ be a simply connected, closed manifold, such that $B\times {\bf R}^{n}$ is
symplectic. Suppose that $f:B\times {\bf R}^{n}\to {\bf R}$  is a smooth function that is
quadratic at infinity. Then $f$ satisfies the conclusion of Corollary \ref{cor:inv} (ii).
If $f$ is a Morse function, then  it satisfies the conclusion of Corollary \ref{cor:inv} (i).

We recall that a function $f$ is quadratic at infinity if
there is a compact set $K\subset B\times {\bf R}^{n}$, a smooth function
$g:B\to {\bf R}$ and a non-degenerate quadratic form $q:{\bf R}^{n}\to {\bf R}$
of index $k$,
$0\leq k\leq n$ such that $f(x,y)=g(x)+q(y)$ for $(x,y)\not\in K$.

The proof of the claim follows from the fact that for a  sufficientely large disk $D\subset
{\bf R}^{n}$ the set $B\times D$ is an isolating neighbourhood for $\gamma(g\oplus q)$ and
it is easy to see that for a big $D$, $(f,B\times D)$ is continuation equivalent to $(g+q,B\times
D)$. Moreover, the flow $\gamma(f)$ is defined for infinite time and this implies that
$S=I_{B\times D}(f)$ verifies property $(\ast)$ by Lemma \ref{lem:prop_ast} (iii).
We also see $d(S)=1$ which is obtained by applying Corollary \ref{cor:cond}.
Indeed, if there is a minimal Morse function continuation equivalent to $(g+q,B\times D^{n})$
one look to the relevant Serre spectral sequence suffices to apply \ref{cor:cond} (v).

Of course, the result is trivial if
$k=0$ or $k=n$ because in these two cases the compactness of $B$  implies
that the level surfaces of $f$ are also compact.

Many examples of symplectic manifolds $B\times {\bf R}^{n}$ are such that
$B$ is itself symplectic and the symplectic form on the product is the product one.
However, there are many examples that are not of this type (one such is the cotangent bundle of
the three sphere) and in all cases the result mentioned is true for an arbitrary symplectic
form on the product.

There are many ways in which to extend this result. We list below a few such possiblities
(that can be also mixed).
\begin{itemize}
\item[(i)] Take $B$ non compact but $g$ such that there is some invariant set for the negative
gradient flow of $g$ that satisfies the conditions (iii) or (iv) in Corollary \ref{cor:cond}
\item[(ii)] Replace $q$ with a function $h:M'\to {\bf R}$ such that the gradient of $h$
induces a flow on $N$, every level hypersurface of $h$ is not compact, the critical
set of $h$ is not trivial and outside of some compact in $N$ the gradient of $h$ is bounded
away from $0$.
\item[(iii)] Replace $B\times {\bf R}^{n}$ with the total space of a fibre bundle
of basis $B$ that is not trivial.
\end{itemize}

Given the result above it is easy to produce non-compact hypersurfaces
$V$ with no compact connected components such that
after a small isotopy they carry a closed characteristic.  Let $M=S^{3}\times {\bf R}^{3}=
\{(x,y,z,w,m,n,p) : x^{2}+y^{2}+z^{2}+w^{2}=1\}$. Take $V\subset M$ be defined
by $V=\{(x,y,z,w,m,n,p) : x+m^{2}+n^{2}-p^{2}=1/2\}$. The symplectic structure on
$M$ comes from the identification $M\approx T^{\ast}S^{3}$.
By a slight variation we may also construct examples of manifolds $M$
that contain connected non-compact hypersurfaces $V_{1}\ldots V_{k}$  pairwise
non homeomorphic and such that after a small  isotopy each contains a closed characteristic.
Whenever Corollary \ref{cor:cond} (v) applies, the number of such hypersurfaces is
(essentially) estimated by half the number of pairs $k,l\in {\bf N}$  which verify the
condition (v) in the corollary. For example, if $M=S^{3}\times S^{3}\times S^{3}\times {\bf
R}^{9}$ (which is symplectic as $M\approx T^{\ast}S^{3}\times T^{\ast}S^{3}\times
T^{\ast}S^{3}$) we can find in $M$ at least two hypersurfaces that are non-compact,
connected and non-homeomorphic that carry a closed characteristic after a small
deformation.
\end{example}

\subsubsection {Final comments and open questions.}

\

A. The one major question that is left open by this paper is to find some
analytic assumptions which, once imposed to an $(f,N)\in\x(M)$ which already
satisfies one of the homotopical conditions in Corollary \ref{cor:cond}, will
imply that $h^{f}$ has some periodic orbits.
One could hope to unify the results on the
existence of periodic orbits which are valid in the compact setting with those that work
in the non-compact case. We should point out however that the major defficiency of the
existence results that we described in the non-compact situation is not so much
their validity only generically but rather the fact that, on one hand, we do not have any control
on the periods of the orbits constructed, and, on the other hand, the genericity
assumption involved is not ``verifiable" (by contrast to the form of genericity
proper to Floer theory where it is determined by the behaviour of a certain
Fredholm operator).

\

B. Of less significance but still of interest are the following questions:

(i) For $(f,N)\in\x(M)$ fixed assume that there is a minimal Morse function
$(g,N')$ continuation equivalent to $(f,N)$
but $\mathcal{C}_{N'}(g)$ has a vanishing differential  (in other words all the relevant
relative Hopf invariants vanish). Find some other homotopical condition on $f$ that
depends on some higher homotopy operation (like Massey products or Toda brackets) in
$\pi_{\ast}^{S}(\Omega M^{+})$ that is sufficient for the existence of a function arbitrarily
close to $g$ whose hamiltonian flow has infinitely many periodic orbits. The first
stage would be to extend Theorem 1 for these higher homotopy operations by relating them
to bordism classes of "maniflods" of possibly broken connecting flow lines. The second
stage would be to extend the argument in Theorem 2 to the case of these bordism classes.

(ii) Extend the arguments in theorems 1, 2  for various other classes of functions,
not necessarily Morse. For example functions having only "neat" singularities as
in \cite{Co3} or Morse-Bott functions.

\

C. Here is a possibly useful extension of Theorem \ref{theo:hopf_conn}
that is not hard to obtain but has not been explicitely included to avoid a considerable
increase in technical difficulty. The statement is similar to that of the theorem but applied to the
Thom map of  $\coprod_{i} Z(P,Q_{i})\hookrightarrow S^{u}(P)$ where
$ind(P)=p>q=ind(Q_{i})$ and $Crit^{N}_{k}(f)=\emptyset$ for $q<k<p$. This Thom map is
equal to the relative Hopf invariant, $H(P,\{Q_{i}\})$, obtained from two succesive
cofibrations $\vee_{i} S^{q-1}_{i}\to N'\to N''$ and the second $S^{p-1}\to N''\to N'''\subset
M$ the first corresponding to the simultaneous cell-attachements corresponding to the
$Q_{i}$'s and the second to the cell attachement corresponding to $P$. Unstably, this Thom
map is richer than the sum of the Thom maps corresponding to the pairs $(P,Q_{i})$ but stably
it equals this sum.   In the statement of Proposition \ref{prop:bded_gen} (i) one can use
$H(P,\{Q_{i}\})$ instead of $H(j_{P},j_{Q})$.

\bibliographystyle{amsplain}

\end{document}